\documentclass[11pt]{amsart}

\usepackage{ulem}
\usepackage{graphicx}
\usepackage{subfigure}
\usepackage{amssymb}
\usepackage[top=25mm, bottom=25mm, left=25mm, right=25mm]{geometry}
\usepackage{color}
\usepackage{bm}
\usepackage{multirow}
\usepackage{epstopdf}
\usepackage{hyperref}
\usepackage{cleveref}

\theoremstyle{definition}
\newtheorem{remark}{Remark}
\newtheorem{definition}{Definition}
\numberwithin{definition}{section}

\newcommand{\pvct}[1]{\bm{#1}}
\newcommand{\vct}[1]{\bm{\mathsf{#1}}}

\newcommand{\pxx}{\pvct{x}}

\newcommand{\uu}{\vct{u}}
\newcommand{\vv}{\vct{v}}

\newcommand{\TT}{\mtx{T}}

\newcommand{\mtx}[1]{\bm{\mathsf{#1}}}

\newcommand{\mtwo}[4]{\left[\begin{array}{cc} #1 & #2 \\ #3 & #4 \end{array}\right]}
\newcommand{\vtwo}[2]{\left[\begin{array}{cc} #1 \\ #2 \end{array}\right]}

\newcommand{\pgnotate}[1]{}

\newcommand{\lsp}{\vspace{3mm}}

\begin{document}

\begin{center}
\textbf{An adaptive high order direct solution technique for elliptic boundary value problems}

\lsp

\textit{\small P.~Geldermans and A.~Gillman}\\
\textit{\small Department of Computational and Applied Mathematics, Rice University}

\lsp

\lsp

\begin{minipage}{143mm}
\textbf{Abstract:} This manuscript presents an adaptive high order
discretization technique for elliptic boundary value problems. 
The technique is applied to an updated version of the Hierarchical Poincar\'e-Steklov 
(HPS) method.  Roughly speaking, the HPS method is based on local pseudospectral 
discretizations glued together with Poincar\'e-Steklov operators.  The 
new version uses a modified tensor product basis which is more efficient and 
stable than previous versions. The adaptive technique exploits the 
tensor product nature of the basis functions to create a criterion
for determining which parts of the domain require additional refinement.
The resulting discretization achieves the user prescribed accuracy and 
comes with an efficient direct solver. The direct solver increases the 
range of applicability to time dependent problems where the cost of 
solving elliptic problems previously limited the use of implicit 
time stepping schemes.
\end{minipage}

\end{center}
%
%
\section{Introduction}

This manuscript presents an adaptive discretization technique for problems
of the form 
\begin{equation}
\label{eq:basic}
\left\{\begin{aligned}
Au(\pxx) =&\ g(\pxx)\qquad &\pxx \in \Omega,\\
u(\pxx)    =&\ f(\pxx)\qquad &\pxx \in \Gamma = \partial \Omega,
\end{aligned}\right.
\end{equation}
where $\Omega$ is a rectangle in $\mathbb{R}^2$ with boundary $\Gamma$,
and where $A$ is a coercive elliptic partial differential operator
\begin{multline}
\label{eq:defA}
[Au](\pxx) = -c_{11}(\pxx)[\partial_{1}^{2}u](\pxx)
-2c_{12}(\pxx)[\partial_{1}\partial_{2}u](\pxx)
-c_{22}(\pxx)[\partial_{2}^{2}u](\pxx)\\
+c_{1}(\pxx)[\partial_{1}u](\pxx)
+c_{2}(\pxx)[\partial_{2}u](\pxx)
+c(\pxx)\,u(\pxx).
\end{multline}

The discretization technique presented here is an 
updated version of the composite spectral discretization techniques 
presented in \cite{2013_martinsson_DtN_linearcomplexity, 2013_martinsson_ItI,2012_spectralcomposite}.  
It is based on local pseudospectral discretizations that are 
``glued'' together by Poincar\'e-Steklov operators.  These Poincar\'e-Steklov
operators are glued in a hierarchically yielding a direct solver.    Hence, the 
discretization technique is called the Hierarchical Poincar\'e-Steklov (HPS) method.  
The adaptive refinement strategy presented in this manuscript is 
inspired by the technique in \cite{lee_greengard_stiff} which 
determines which parts of the geometry to refine by looking at 
Chebychev expansion coefficients of the local approximate solution.  
Like the HPS methods in \cite{2013_martinsson_DtN_linearcomplexity, 
2013_martinsson_ItI,2012_spectralcomposite,2016_bodyload},
the adaptive discretization technique can also be modified to handle a 
range of different domains, including
curved ones.  Additional novelty of this 
paper lies in an update to the local discretization. The new local discretization
uses a \emph{modified} tensor product basis which makes the local discretization 
less expensive than previous versions \cite{2013_martinsson_DtN_linearcomplexity, 
2013_martinsson_ItI,2012_spectralcomposite,2016_bodyload} and the whole algorithm 
easier to implement.

While constructing the adaptive discretization and the direct solver has 
a computational cost that scales $O(N^{3/2})$ where $N$ is the number of 
discretization points, the cost of applying the solver is $O(N\log N)$ with 
a small constant.  The constant in the solve step is typically much 
smaller than for a uniform discretization thus making the method useful 
for applications that involve \textit{many} elliptic solves that require locally 
refined high order discretizations.  For example, having an 
efficient direct solver for elliptic partial differential equations can increase
the range of problems for which implicit time stepping schemes are 
computationally affordable.

\subsection{Overview of discretization technique}
Roughly speaking the adaptive discretization technique can be broken into three steps.
\begin{itemize}
 \item[Step 1:] First, the geometry is partitioned into a collection 
 of patches using an quad tree
 with an adaptive interpolation strategy applied so that the coefficients in 
 (\ref{eq:defA}) and body load function $g(\pxx)$ in (\ref{eq:basic}) are captured to the
 user prescribed tolerance $\epsilon$.
 \item[Step 2:] Next each patch is discretized using a high order spectral 
 collocation technique and the patches are ``glued'' together at
 the boundaries via a Poincar\'e-Steklov operators in a hierarchical
 fashion.  In the process of gluing patches together, solution operators that 
 propagate boundary data to the interior of a box are constructed.  
 Then by applying the solution operators (small matrix vector multiplies) the 
 boundary data is propagated down the hierarchical 
 tree giving an approximate solution on each patch.
 \item[Step 3:] All patches are checked to see if they need to be further 
 refined.  If there are patches marked for refinement, they are 
 refined and steps 2 and 3 are repeated until no patches are 
 marked for refinement.  If the refinement is localized in the domain, the 
 bulk of the computation from step 2 can be reused.
\end{itemize}

While the method can be employed with any Poincar\'e-Steklov operator, for simplicity 
of presentation, this paper uses the Dirichlet-to-Neumann operator for gluing boxes as 
in \cite{2013_martinsson_DtN_linearcomplexity, 2012_spectralcomposite, 2016_bodyload}.  
For the Helmholtz experiments in this paper, the impedance-to-impedance (ItI) operator 
is used instead.  \cite{2013_martinsson_ItI} presents the ItI version of the solution 
technique for a homogeneous PDE. The appendix of this manuscript presents the ItI based solution technique when
there is a body load $g(\pxx)$.  

\subsection{Applications utilizing the HPS method}
While the HPS method is relatively new, it is already being utilized for scattering 
problems.  Applications involving scattering problems include underwater acoustics \cite{underwater},
ultrasound and microwave tomography \cite{3dmicrowave,wadbro10},
wave propagation in metamaterials and photonic crystals, and seismology \cite{seismic3D}.
In \cite{2013_martinsson_ItI}, the HPS method was extended to free space scattering problems
where the deviation from a constant coefficient problem had compact support. 
The numerical results in that paper showed the method did not observe pollution for 
problems where the support of the deviation from constant coefficient was 100 time the 
smallest wavelength in size.  In \cite{2017_Borges}, the method was utilized to build a inverse scattering 
solver via the recursive linearization procedure proposed in \cite{1995_chen}. The 
recursive linearization procedure requires solving a sequence of linear least squares 
problems at successively higher frequencies to reconstruct an unknown sound speed.  Next, 
in \cite{2018_Borges}, the HPS method was utilized for inverse scattering problems 
with a random noisy background medium.  In each of these inverse scattering solvers,
the least squares solve requires solving the same variable coefficient elliptic partial 
differential equation many times to apply the forward and adjoint operators. 
The proposed adaptive discretization could improve the efficiency of the techniques listed 
in this section.

\subsection{Prior and related work}
There is a vast literature of adaptive methods for finite element (FEM) based 
discretization techniques for elliptic problems.  A high level overview is 
presented here.  Early works \cite{1988_Eriksson,1978_bab} focused on 
defining appropriate error estimators for Poisson problems using face and volume residuals 
giving the user the ability to identify where to refine.  Recent trends in adaptive 
FEM for Poisson problems focus on proving that the adaptive algorithms converge \cite{2002_Morin,2006_AMFEM}.
A local indicators and error estimators for FEM applied to Helmholtz problems are presented in \cite{1997_bab}.
There has also been an extensive work on hp-adaptivity \cite{2000_rachowicz,dem_book1,dem_book2}.
The adaptive discretization presented here is an an h-adaptive scheme which is specific for the HPS 
discretization technique.  The local error indicator can be (and is) applied to both Poisson and 
Helmholtz problems.  The relative convergence error stopping criterion determines if the problem
has been resolved.

The direct solver for the HPS discretization is related to the direct solvers for sparse systems
arising from finite difference and finite element discretizations of elliptic
PDEs such as the classical nested dissection method of George \cite{george_1973,hoffman_1973}
and the multifrontal methods by Duff and others \cite{1989_directbook_duff}.  These 
methods can be viewed as a hierarchical version of the
``static condensation'' idea in finite element analysis \cite{1974_wilson_static_condensation}.
High order finite difference and finite element discretizations lead to large frontal matrices
(since the ``dividers'' that partition the grid have to be wide), and consequently
very high cost of the LU-factorization (see, e.g., Table 2 in \cite{2013_martinsson_DtN_linearcomplexity}).  
It has been demonstrated that the dense matrices that arise in these solvers
have internal structure that allows the direct solver to be accelerated to linear or close to linear
complexity, see, e.g., \cite{2009_xia_superfast,Adiss,2007_leborne_HLU,2009_martinsson_FEM,2011_ying_nested_dissection_2D}.
The HPS discretization technique has one dimensional ``dividers'' independent of order and thus
the direct solver only pays (in terms of computational complexity) the price of high order 
at the lowest level in the hierarchical tree.  
The same ideas that accelerate the nested dissection and multifrontal solvers can be applied 
the HPS direct solver \cite{2013_martinsson_DtN_linearcomplexity}.

In the previous versions of the HPS method special care was taken to deal with or avoid
discretization points at the corners of the small patches.  The method presented
in \cite{2012_spectralcomposite} involves tedious bookkeeping of corner points.
Additionally, possible singularities at the corners of the geometry $\Omega$ are
of concern.  By introducing interpolation at the level of the local discretizations,
the methods in \cite{2013_martinsson_DtN_linearcomplexity,2016_bodyload, 2013_martinsson_ItI}
avoid the corners of $\Omega$.  The new local discretization presented in this 
manuscript does not involve the corner points at all; thus improving the 
robustness and efficiency of the method.

\subsection{Outline of paper}
For simplicity of presentation, the proposed algorithm is described for a PDE with 
no body load (i.e. $g(\pxx) = 0$ in (\ref{eq:basic})).  The manuscript begins 
by reviewing the HPS
method with uniform refinement in section \ref{sec:HPS} but with the 
\textit{new} local discretization technique.  
Next, the adaptive refinement procedure
is presented in section \ref{sec:adaptive}.  Then, numerical experiments 
demonstrating the performance of the method in section \ref{sec:numerics}.
Finally the manuscript concludes with a summary of the paper in section 
\ref{sec:summary}.

\section{The HPS method}
\label{sec:HPS} 
This section presents the HPS method with a new local discretization technique. 
The HPS method begins by partitioning the domain $\Omega$ into a 
collection of square (or possibly rectangular) boxes, called \textit{leaf boxes}.
Throughout this paper, we assume that the parameter for the order of the 
discretization $n_c$ is fixed ($n_c=16$ is often a good choice).  
For a uniform discretization, the size of all leaf boxes is chosen so that 
any potential $u$ of equation (\ref{eq:basic}), as well as its first 
and second derivatives, can be accurately interpolated from their values 
at the local discretization points on any leaf box.

\begin{figure}

\setlength{\unitlength}{1mm}
\begin{picture}(120,40)
\put(05,00){
\includegraphics[width=120mm]{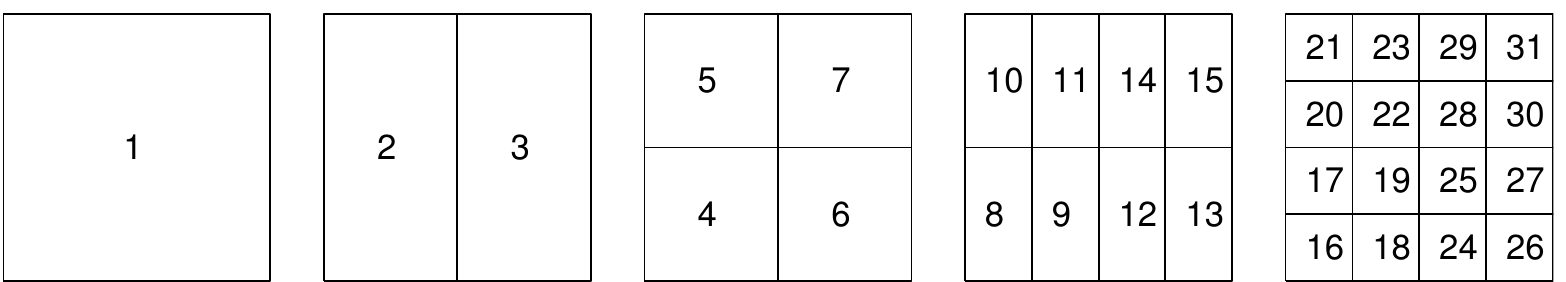}}
\end{picture}
\caption{\label{fig:tree_numbering}
The square domain $\Omega$ is split into $4 \times 4$ leaf boxes.
These are then gathered into a binary tree of successively larger boxes
as described in Section \ref{sec:HPS}. One possible enumeration
of the boxes in the tree is shown, but note that the only restriction is
that if box $\tau$ is the parent of box $\sigma$, then $\tau < \sigma$.}
\end{figure}

Next a binary tree on the collection of leaf boxes is constructed by
hierarchically merging them, making sure that all boxes on
the same level are roughly of the same size, cf.~Figure
\ref{fig:tree_numbering}.  The boxes should be ordered so
that if $\tau$ is a parent of a box $\sigma$, then $\tau < \sigma$. We
also assume that the root of the tree (i.e.~the full box $\Omega$) has
index $\tau=1$. We let $\Omega^{\tau}$ denote the domain associated with box $\tau$.
If a box $\rho$ is child of $\sigma$ and $\sigma$ is a child of $\tau$, 
we call $\rho$ a \textit{grandchild} of $\tau$.  For example in Figure 
\ref{fig:tree_numbering}, boxes $16-19$ are grandchildren of box $4$.  
(This vocabulary is needed for the adaptive scheme presented in section \ref{sec:adaptive}.)

For each leaf box, approximate Dirichlet-to-Neumann (DtN) and 
solution operators are constructed via the \textit{modified} spectral
collocation method presented in section \ref{sec:leaf}.  The DtN approximations
are ``glued'' together in a hierarchical fashion two boxes at a time. 
Section \ref{sec:merge} presents the technique for constructing approximate DtN
and solution operators for the union of two boxes.  Algorithm 1
gives an overview of the construction of the discretization and direct solver.
Once the hierarchical collection of approximate solution operators 
is constructed, the solution on the interior can be found for $O(N\log N)$ cost
via Algorithm 2.

\begin{definition}[Dirichlet-to-Neumann map]
For domain $\Omega$ with boundary $ \Gamma$, the Dirichlet-to-Neumann (DtN) operator  
$T:H^1(\Gamma)\to L^2(\Gamma)$ is defined by
\begin{equation}
T f\ = u_n,
\label{def:T}
\end{equation}
for any Dirichlet boundary data $f(x)\in H^1(\Gamma)$, where 
$u_n$ denotes the normal derivative of $u$ on $\Gamma$ in the direction of the normal vector $n$ pointing
out of $\Omega$.
\label{def:dtn}
\end{definition}


\subsection{Leaf computation}
\label{sec:leaf}

This section describes a modified spectral collocation method for 
constructing approximate DtN $\mtx{T}^\tau$ and solution $\mtx{\Psi}^\tau$ 
operators for a leaf box $\tau$.  

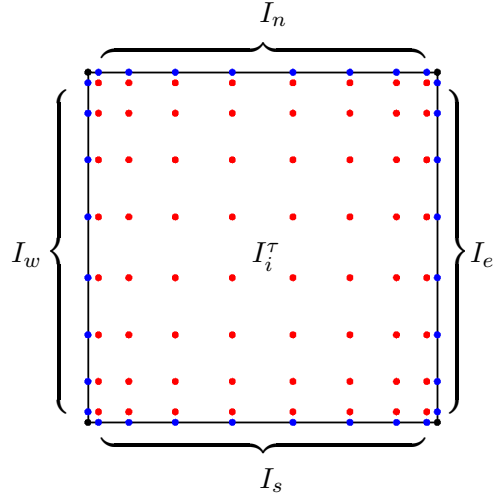
\begin{figure} 
\centering
\setlength{\unitlength}{1mm}
\begin{picture}(60,60)
\put(05,05){\includegraphics[height = 50mm]{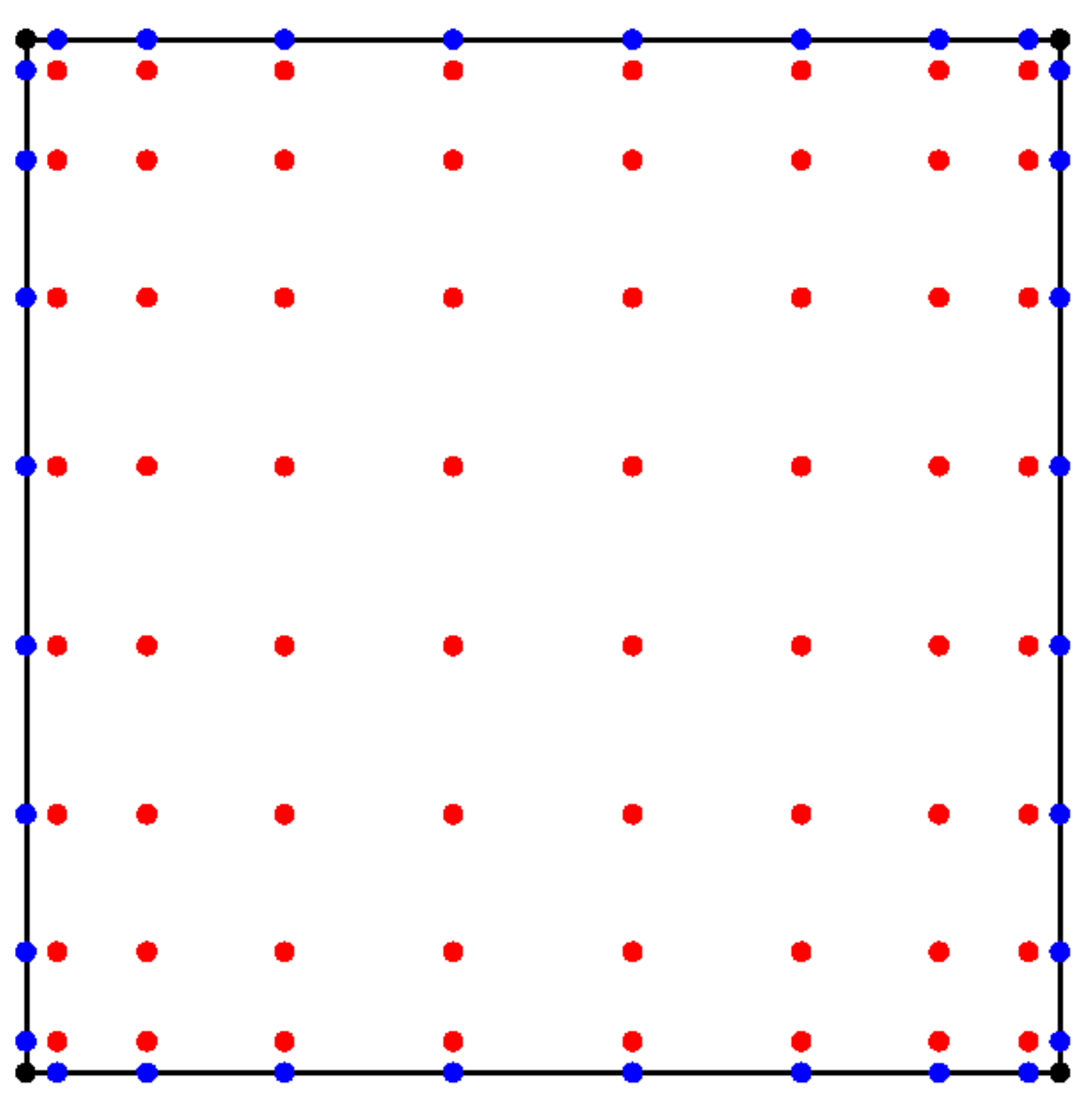}}
\put(28,28){$I^\tau_i$}
\put(8,05){$\underbrace{\hspace{4.3cm}}$}
\put(29,-2){$I_s$}
\put(08,55){$\overbrace{\hspace{4.3cm}}$}
\put(29,60){$I_n$}
\put(54,08){{\rotatebox{90}{$\underbrace{\hspace{4.3cm}}$}}}
\put(57,28){$I_e$}
\put(1,08){{\rotatebox{90}{$\overbrace{\hspace{4.3cm}}$}}}
\put(-4,28){$I_w$}
\end{picture}
\caption{\label{fig:leaf} Illustration of the discretization 
points for a leaf box $\tau$.  The points in blue are the 
boundary points with indices $I^\tau_b = [I_s,I_e,I_n,I_w]$.
The points in red are the interior points with indices $I^\tau_i$.  
The points in black are the omitted corner points.
}
\end{figure}

The modified spectral collocation technique begins with the 
classic $n_c\times n_c$ product Chebychev grid and the 
corresponding differential matrices $\mtx{D}_x$ and $\mtx{D}_y$ 
from for example \cite{2000_trefethen_spectral_matlab}.   	
Let ${I}^\tau_i$ denote the index vector corresponding to points on the 
interior of $\Omega^\tau$ and ${I}^\tau_b$ denote the index vector 
corresponding to points on the boundary of $\Omega^\tau$ \textbf{not} 
including the corner points based on the tensor classic tensor grid.
Figure \ref{fig:leaf} illustrates the 
indexing of the points in terms of the classic discretization.  Thus 
$\{\pxx_j\}_{j=1}^{n_c^2-4}$ denotes the discretization points 
in $\Omega^\tau$ given by the union of the red and blue points 
in Figure \ref{fig:leaf}.  
We order the solution vector $\vct{u}$ and flux
vector $\vct{v}$ according to the following:
$\vct{u} = \vtwo{\vct{u}_b}{\vct{u}_i}$ where $\vct{u}_b$ and $\vct{u}_i$ 
denote the approximate values of the solution on the 
boundary and the interior, respectively.  The ordering of the 
entries related to the boundary corresponding to the discretization points 
is $I^\tau_b = [I_s,I_e,I_n,I_w]$ where $I_s$ denotes the blue points on the south
boundary in Figure \ref{fig:leaf}, etc. Let $I^\tau = [I^\tau_b,I^\tau_i]$ denote
the collection of all indices that are used in the discretization.

Thanks to the tensor product basis, we know the entries of $\mtx{D}_x$ and 
$\mtx{D}_y$ corresponding to the interaction of the corner points with 
the points on the interior of $\Omega^\tau$ are zero.  The directional
basis functions for the other points on the boundary are not impacted by 
the removal of the corner points.  Thus the differential operators from
the classic pseudospectral discretization can be used to create the
approximation of the local differential operator and DtN.

The classic discrete approximation of the differential operator on $\Omega^\tau$ is given by 
\begin{equation*}
\mtx{A} =
-\mtx{C}_{11}\mtx{D}_x^{2}
-2\mtx{C}_{12}\mtx{D}_x\mtx{D}_y
-\mtx{C}_{22}\mtx{D}_y^{2}
+\mtx{C}_{1}\mtx{D}_x
+\mtx{C}_{2}\mtx{D}_y
+\mtx{C},
\end{equation*}
where $\mtx{C}_{11}$ is the diagonal matrix with diagonal entries $\{c_{11}(\pxx_{k})\}_{k=1}^{n_c^{2}}$,
and the other matrices $\mtx{C}_{ij}$, $\mtx{C}_{i}$, $\mtx{C}$ are defined analogously.
Then the discretized differential equation on the new set of discretization 
points is given by 
$$ \left[\begin{array}{c}\mtx{I}\quad \mtx{0} \\  \mtx{A}_{i,b} \quad \mtx{A}_{i,i} \end{array}
\right]
\vtwo{\vct{u}_b}{\vct{u}_i}= \vtwo{\hat{\vct{f}}}{\vct{0}}$$
where $\mtx{A}_{i,i} = \mtx{A}({I}^\tau_i,{I}^\tau_i)$ is a matrix of size
${(n_c-2)^2\times (n_c-2)^2}$, 
$\mtx{A}_{i,b} = \mtx{A}({I}^\tau_i,{I}^\tau_b)$ is a matrix of size $(n_c-2)^2 \times (4n_c-8)$, 
and $\hat{\vct{f}}$ is vector of length $4n_c-8$ containing fictitious Dirichlet boundary data.

When the boundary data is known, the approximate solution at the interior points is given by
\begin{equation}
\vct{u}_i = -\mtx{A}_{i,i}^{-1}\mtx{A}_{i,b}\vct{u}_b = \mtx{\Psi}^\tau\vct{u}_b
\label{eq:soln} 
\end{equation}
where the matrix $\mtx{\Psi}^\tau$ is the approximate solution operator. Since the 
matrix $\mtx{A}_{i,i}$ is not large (even for $n_c = 16$), 
it can be inverted quickly using dense linear algebra.

Let $\mtx{L}$ denote the matrix made up of four block row matrices 
corresponding to taking the normal derivative of the basis functions
on the leaf $\tau$ along each of the edges.  In terms of the discrete
operators $\mtx{L}$ is given by 
$$\mtx{L} = \left[\begin{array}{c} \mtx{D}_x(I_s,I^\tau) \\
                           \mtx{D}_y(I_e,I^\tau)\\
                           \mtx{D}_x(I_n,I^\tau)\\
                           \mtx{D}_y(I_w,I^\tau)\end{array}\right].$$
To construct the approximate DtN operator $\mtx{T}^\tau$, we take 
the normal derivative of the solution by applying $\mtx{L}$ to 
$\vtwo{\mtx{I}_{4n_c-8}}{\mtx{\Psi}^\tau}$, i.e.
 $$\mtx{T}^\tau = \mtx{L}\vtwo{\mtx{I}_{4n_c-8}}{\mtx{\Psi}^\tau}$$
where $\mtx{I}_{4n_c-8}$ denotes the identity matrix of size $4n_c-8$. 

\begin{remark}
 The classic tensor product discretization can be used to formulate the 
 new discretization thanks to the separable basis (i.e. the corner points
 do not contribute the discretized differential equation).  
 While interpolation along the edges without the corners is less accurate 
 than if the corners were included, it is stable \cite{2006_smith}.  
 Since the discretization is run at high order (typically $n_c \geq 16$), 
 a loss in accuracy is not observed in practice. 
\end{remark}

\subsection{Merging two boxes}
\label{sec:merge}
This section reviews of the procedure for constructing the DtN and solution 
matrices for the union of two boxes for which DtN matrices have already 
been constructed. More detailed descriptions are presented in 
\cite{2013_martinsson_DtN_linearcomplexity, 2012_spectralcomposite,
2016_bodyload}.

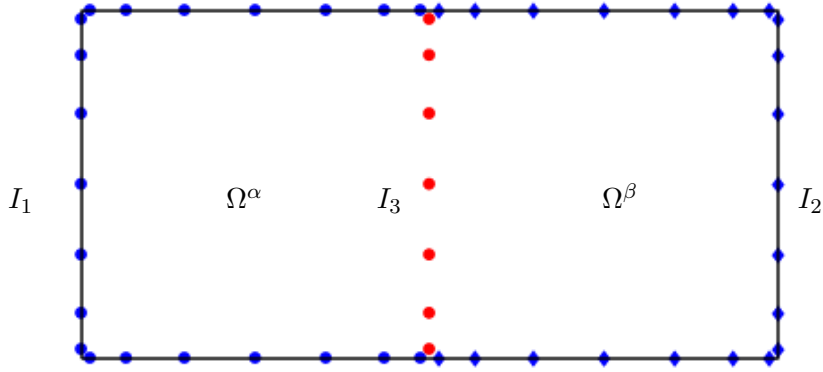
\begin{figure}
\centering
\setlength{\unitlength}{1mm}
\begin{picture}(95,55)
 \put(00,00){\includegraphics[height=55mm]{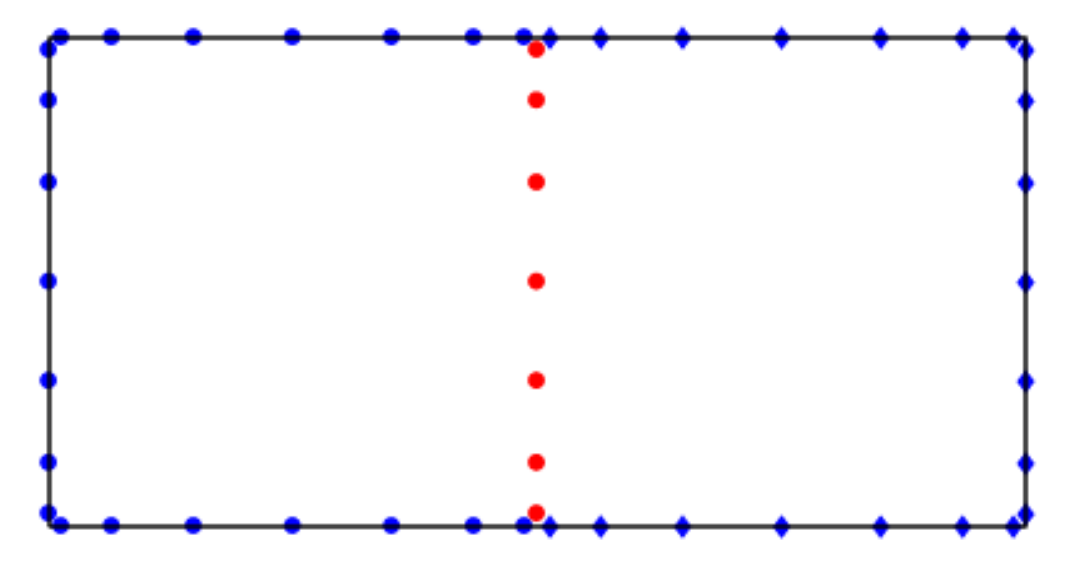}}
\put(24,25){$\Omega^{\alpha}$}
\put(74,25){$\Omega^{\beta}$}
\put(-5,25){$I_{1}$}
\put(100,25){$I_{2}$}
\put(44,25){$I_{3}$}
\end{picture}
\caption{Notation for the merge operation described in Section \ref{sec:merge}.
The rectangular domain $\Omega$ is formed by two squares $\Omega^{\alpha}$ 
and $\Omega^{\beta}$.  The sets $I_{1}$ (blue circles) and $I_{2}$ (blue diamonds) 
form the exterior nodes, while $I_{3}$ (red circles) consists of the interior nodes.}
\label{fig:siblings_notation}
\end{figure}

Let $\Omega^{\tau}$ denote a box with children $\Omega^{\alpha}$ and
$\Omega^{\beta}$ so that
$$
\Omega^{\tau} = \Omega^{\alpha} \cup \Omega^{\beta}.$$ 
For concreteness, let us assume that $\Omega^{\alpha}$ and
$\Omega^{\beta}$ share a vertical edge as shown in Figure \ref{fig:siblings_notation}.  
We partition the points on $\partial\Omega^{\alpha}$ and $\partial\Omega^{\beta}$ into three sets:
\begin{tabbing}
\mbox{}\hspace{5mm}\= $I_{1}$ \hspace{4mm} \=
Boundary nodes of $\Omega^{\alpha}$ that are not boundary nodes of $\Omega^{\beta}$.\\
\> $I_{2}$ \> Boundary nodes of $\Omega^{\beta}$ that are not boundary nodes of $\Omega^{\alpha}$.\\
\> $I_{3}$ \> Boundary nodes of both $\Omega^{\alpha}$ and $\Omega^{\beta}$ that are \textit{not} boundary nodes of the\\
\> \>  union box $\Omega^{\tau}$.
\end{tabbing}
The indexing for the points on the interior and boundary of 
$\Omega^\tau$ are $I_{\rm i}^\tau = I_3$ and $I_{\rm b}^\tau = [I_1,I_2]$, 
respectively.

Let $u$ denote a solution to (\ref{eq:basic}), with tabulated potential values $\uu$
and boundary fluxes $\vv$. Ordering the DtN operators according to the $I_k$ defined
 in Figure \ref{fig:siblings_notation} results in the equations
\begin{equation}
\label{eq:bittersweet1}
\left[\begin{array}{c}
\vv_{1}\\ \vv_{3}
\end{array}\right] =
\left[\begin{array}{ccc}
\mtx{T}_{1,1}^{\alpha} & \mtx{T}_{1,3}^{\alpha} \\
\mtx{T}_{3,1}^{\alpha} & \mtx{T}_{3,3}^{\alpha}
\end{array}\right]\,
\left[\begin{array}{c}
\uu_{1}\\ \uu_{3}
\end{array}\right],
\qquad {\rm and} \qquad
\left[\begin{array}{c}
\vv_{2}\\ \vv_{3}
\end{array}\right] =
\left[\begin{array}{ccc}
\mtx{T}_{2,2}^{\beta} & \mtx{T}_{2,3}^{\beta} \\
\mtx{T}_{3,2}^{\beta} & \mtx{T}_{3,3}^{\beta}
\end{array}\right]\,
\left[\begin{array}{c}
\uu_{2}\\ \uu_{3}
\end{array}\right]
\end{equation}
where $\mtx{T}^\alpha_{1,1} = \mtx{T}^\alpha(I_1,I_1)$, etc. 
Noting that $\vv_3$ and the solution $\uu_3$ is the 
same for each box (since the solution is smooth), the solution operator
$\mtx{\Psi}^\tau$ is found by equating the bottom two row equations
of (\ref{eq:bittersweet1});
\begin{equation}
\uu_{3} =\ \bigl(\TT^{\alpha}_{3,3} - \TT^{\beta}_{3,3}\bigr)^{-1}
\bigl[-\TT^{\alpha}_{3,1}\ \big|\ \TT^{\beta}_{3,2}] \left[\begin{array}{c}\uu_{1} \\ \uu_{2}\end{array}\right]  
=\
\mtx{\Psi}^{\tau}\,
\left[\begin{array}{c}\uu_{1} \\ \uu_{2}\end{array}\right].
 \label{eq:Psi}
\end{equation}
The operator $\mtx{\Psi}^\tau$ in (\ref{eq:Psi}) maps the solution on the boundary of
$\Omega^\tau$ to the solution  $\mtx{u}_3$ on the interior edge.  
The DtN operator is then constructed by plugging equation (\ref{eq:Psi}) into
the top row equations in (\ref{eq:bittersweet1}) and combining them to a 
matrix equation.  The result is 
\begin{align}
\label{eq:DtN}
\left[\begin{array}{c}\vv_{1} \\ \vv_{2}\end{array}\right] = \mtx{T}^{\tau}\,
\left[\begin{array}{c}\uu_{1} \\ \uu_{2}\end{array}\right]
\end{align}
where 
\begin{equation}
\label{eq:merge}
\mtx{T}^{\tau} =
\left[\begin{array}{ccc}
\TT_{1,1}^{\alpha} & \mtx{0} \\
\mtx{0} & \TT_{2,2}^{\beta }
\end{array}\right] +
\left[\begin{array}{c}
\TT_{1,3}^{\alpha} \\
\TT_{2,3}^{\beta}
\end{array}\right]\,
\bigl(\TT^{\alpha}_{3,3} - \TT^{\beta}_{3,3}\bigr)^{-1}
\bigl[-\TT^{\alpha}_{3,1}\ \big|\ \TT^{\beta}_{3,2}\bigr].
\end{equation}

\begin{figure}[h!]
\fbox{
\begin{minipage}{120mm}
\begin{center}
\textsc{Algorithm 1} (build solution operators)
\end{center}

This algorithm builds the global Dirichlet-to-Neumann operator for (\ref{eq:basic}).\\
It also builds all the solution matrices $\mtx{\Psi}^{\tau}$ required for 
constructing an approximation to $u$ at any interior point.\\
It is assumed that if node $\tau$ is a parent of node $\sigma$, then $\tau < \sigma$.

\rule{\textwidth}{0.5pt}

\begin{tabbing}
\mbox{}\hspace{7mm} \= \mbox{}\hspace{6mm} \= \mbox{}\hspace{6mm} \= \mbox{}\hspace{6mm} \= \mbox{}\hspace{6mm} \= \kill
(1)\> \textbf{for} $\tau = N_{\rm boxes},\,N_{\rm boxes}-1,\,N_{\rm boxes}-2,\,\dots,\,1$\\
(2)\> \> \textbf{if} ($\tau$ is a leaf)\\
(3)\> \> \> Construct $\mtx{T}^{\tau}$ and $\mtx{\Psi}^\tau$ via the process described in Section \ref{sec:leaf}.\\
(4)\> \> \textbf{else}\\
(5)\> \> \> Let $\sigma_{1}$ and $\sigma_{2}$ be the children of $\tau$.\\
(6)\> \> \> Split $I_{\rm b}^{\sigma_{1}}$ and $I_{\rm b}^{\sigma_{2}}$ into vectors $I_{1}$, $I_{2}$, and $I_{3}$ as shown in Figure \ref{fig:siblings_notation}.\\
(7)\> \> \> $\mtx{\Psi}^{\tau} = \bigl(\TT^{\sigma_{1}}_{3,3} - \TT^{\sigma_{2}}_{3,3}\bigr)^{-1}
                           \bigl[-\TT^{\sigma_{1}}_{3,1}\  \big|\
                                  \TT^{\sigma_{2}}_{3,2}\bigr]$\\
(8)\> \> \> $\TT^{\tau} = \left[\begin{array}{ccc}
                          \mtx{T}_{1,1}^{\sigma_{1}} & \mtx{0}\\
                          \mtx{0} & \mtx{T}_{2,2}^{\sigma_{2} }
                          \end{array}\right] +
                    \left[\begin{array}{c}
                          \TT_{1,3}^{\sigma_{1}} \\
                          \TT_{2,3}^{\sigma_{2}}
                          \end{array}\right]\,\mtx{\Psi}^{\tau}$.\\
(9)\> \> \> Delete $\TT^{\sigma_{1}}$ and $\TT^{\sigma_{2}}$.\\
(10)\> \> \textbf{end if}\\
(11)\> \textbf{end for}
\end{tabbing}
\end{minipage}}
\label{alg:precompute}
\end{figure}

\begin{figure}[h!]
\fbox{
\begin{minipage}{120mm}
\begin{center}
\textsc{Algorithm 2} (solve BVP once solution operator has been built)
\end{center}

This program constructs an approximation $\uu$ to the solution $u$ of (\ref{eq:basic}).\\
It assumes that all matrices $\mtx{\Psi}^{\tau}$ have already been constructed in a pre-computation.

\rule{\textwidth}{0.5pt}

\begin{tabbing}
\mbox{}\hspace{7mm} \= \mbox{}\hspace{6mm} \= \mbox{}\hspace{6mm} \= \mbox{}\hspace{6mm} \= \mbox{}\hspace{6mm} \= \kill
(1)\> $\uu(k) = f(\pxx_{k})$ for all $k \in I_{\rm b}^{1}$.\\
(2)\> \textbf{for} $\tau = 1,\,2,\,3,\,\dots,\,N_{\rm boxes}$\\
(3) \> \> $\uu(I_{\rm i}^{\tau}) = \mtx{\Psi}^{\tau}\,\uu(I_{\rm b}^{\tau})$.\\
(4)\> \textbf{end for}
\end{tabbing}

\vspace{1mm}
\end{minipage}}
\label{alg:Solve}
\end{figure}

\section{Adaptive discretization}
\label{sec:adaptive}
This section presents an adaptive discretization technique for the boundary value 
problem (\ref{eq:basic}) where the coefficient functions, right hand side and
boundary data are smooth functions.  As with the uniform discretization 
technique, the adaptive method produces a direct solver.  The approximate solution 
obtained from the adaptive procedure is accurate  
(in the relative error) to a user prescribed tolerance $\epsilon$.

At a high level, the idea stems from the fact that the discretization on a leaf can 
be accurate enough to capture the solution locally if it was given correct
boundary data.  The indicator for further refinement we propose in this section 
determines if the local basis is good enough to capture the solution locally.  
The stopping criterion for the 
adaptive procedure is based on the relative convergence error.  This ensures that
each leaf is given accurate boundary data.  Starting such an adaptive technique
with a global discretization of $\Omega$ would be computationally prohibitive.
Instead, we initialize the mesh by utilizing the fact that the basis on a leaf should be able to 
represent the coefficient functions and the right hand side in (\ref{eq:basic})
to the user prescribed $\epsilon$.  

\begin{remark}
 In practice, one could likely get away with asking for less accuracy of the 
 adaptive interpolation scheme. Since the interpolation is inexpensive
 compared to the cost of building the discretization and direct solver, we 
 choose to be cautious.
\end{remark}

The algorithm can be broken into seven steps.

\begin{itemize}
 \item[Step 1:] Use the adaptive interpolation technique from section \ref{sec:interp} 
 applied to the coefficient functions and the right hand side of equation (\ref{eq:basic}).
 This yields an initial mesh.
 \item[Step 2:] Construct an HPS solver for the non-uniform mesh resulting from the 
 adaptive interpolation scheme via the techniques presented in section \ref{sec:nonuniform}.\\ 
 Let $\vct{u}^{\tau,{\rm old}}$  denote the approximate solution on leaf box $\tau$.  
 \item[Step 3:] Use indicator presented in section \ref{sec:mark} to determine which 
 boxes need additional refinement. 
 \item[Step 4:]  If a leaf box $\tau$ has been 
marked for refinement, split into into four boxes ($\alpha_1$, $\alpha_2$, $\alpha_3$ and $\alpha_4$). 
\item[Step 5:] Discretize the new leaf boxes and update the direct solver.  
Since the discretization is 
localized and the direct solver is naturally domain decomposing, the direct solver
can efficiently be updated without touching the entire geometry (see section \ref{sec:update}).\\
Let $\vct{u}^{\tau,{\rm new}}$ denote the solution on leaf box $\tau$ obtained with the new mesh. 
\item[Step 6:] Check the relative convergence error by sweeping over all
the leaf boxes on the old tree.\\
If leaf box $\tau$ was not refined, the relative convergence 
error for that box is defined to be 
$$E^\tau_{\rm rel} = \frac{\|\vct{u}^{\tau,{\rm old}}-\vct{u}^{\tau,{\rm new}}\|_2}{\|\vct{u}^{\tau,{\rm old}}+\vct{u}^{\tau,{\rm new}}\|_2}.$$ 
If leaf box $\tau$ was refined with grandchildren $\alpha_1$, $\alpha_2$, $\alpha_3$ and $\alpha_4$
then the relative convergence error is defined as
$$E^\tau_{\rm rel} = \frac{\|\vct{u}^{\tau,{\rm old}}-\mtx{L}\vct{u}^{\tau}_{\rm fine}\|_2}{\|\vct{u}^{\tau,{\rm old}}+\mtx{L}\vct{u}^{\tau}_{\rm fine}\|_2}$$
where $\vct{u}^\tau_{\rm fine} = \left[\begin{array}{c}{\vct{u}^{\alpha_1,{\rm new}}}\\{\vct{u}^{\alpha_2,{\rm new}}}\\
                                        {\vct{u}^{\alpha_3,{\rm new}}}\\ {\vct{u}^{\alpha_4,{\rm new}}} 
                                        \end{array}\right]$, and 
 $\mtx{L}$ is a matrix that interpolates functions from the fine discretization points to the coarse 
 discretization points.  \\
\item[Step 7:] If the average relative $E_{\rm rel} = \frac{1}{n_\tau}\sum_{\tau}E^\tau_{\rm rel}> \epsilon$ where $\tau$ 
 is a leaf box on the old tree and $n_\tau$ is the number of leaf boxes in the old tree, the algorithm terminates.  
 Otherwise, the vectors $\vct{u}^{\tau,{\rm new}}$ get 
 the label ${\bf old}$ and return to Step 3.

\end{itemize}

\subsection{Adaptive interpolation}
\label{sec:interp} 
In order to keep the cost of the adaptive discretization as low as possible, we 
first create a mesh which allows for the smooth functions in (\ref{eq:basic}) 
to be approximated with the local bases to the user prescribed tolerance
$\epsilon$.  For simplicity of presentation, we describe the technique
for interpolating a general smooth function $f(\pxx)$ on $\Omega$.
  
\begin{figure}[h!] 
\centering
\setlength{\unitlength}{1mm}
\begin{picture}(95,55)
\put(-10,00){\includegraphics[height = 50mm]{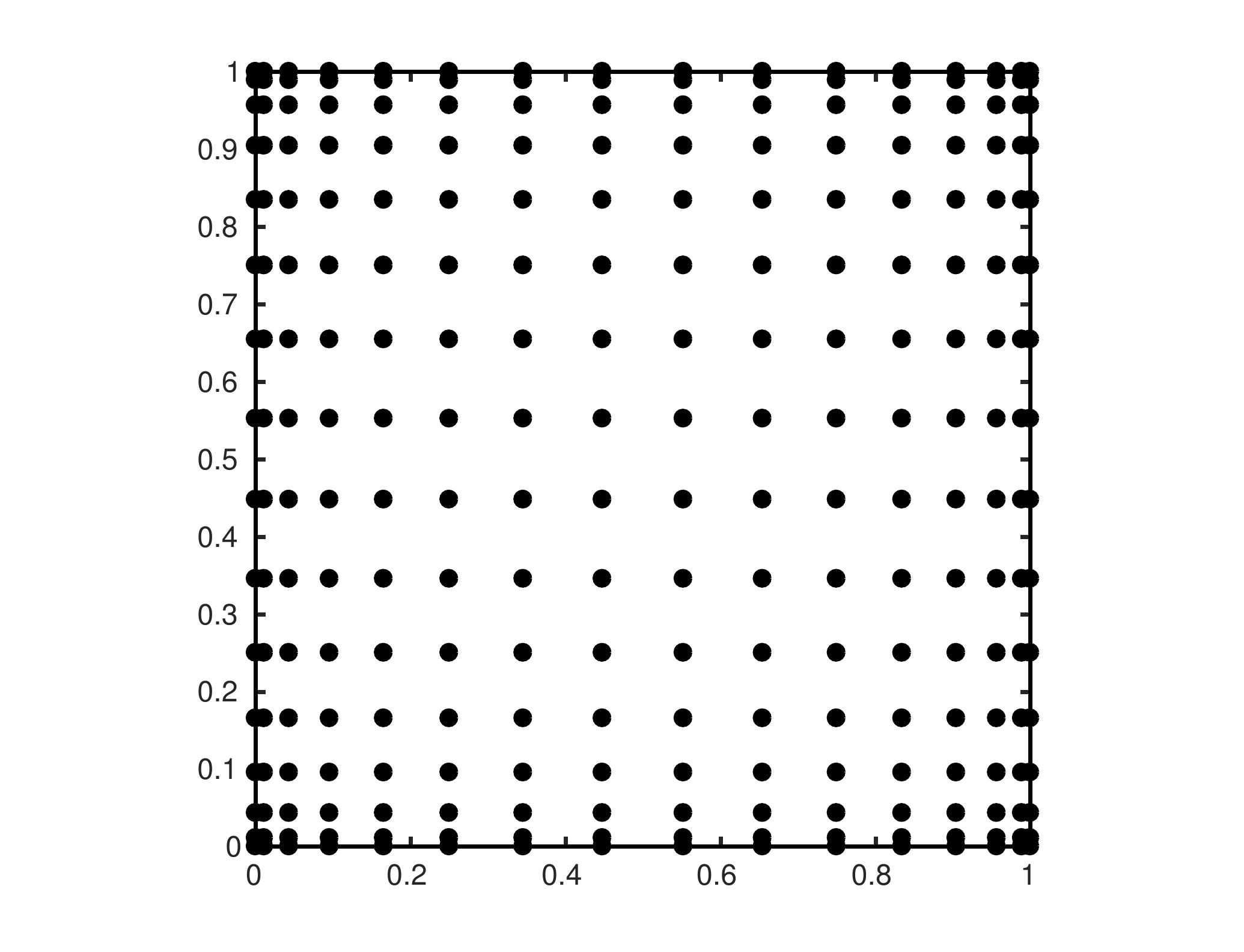}}
\put(21,-3){(a)}
\put(50,00){\includegraphics[height = 50mm]{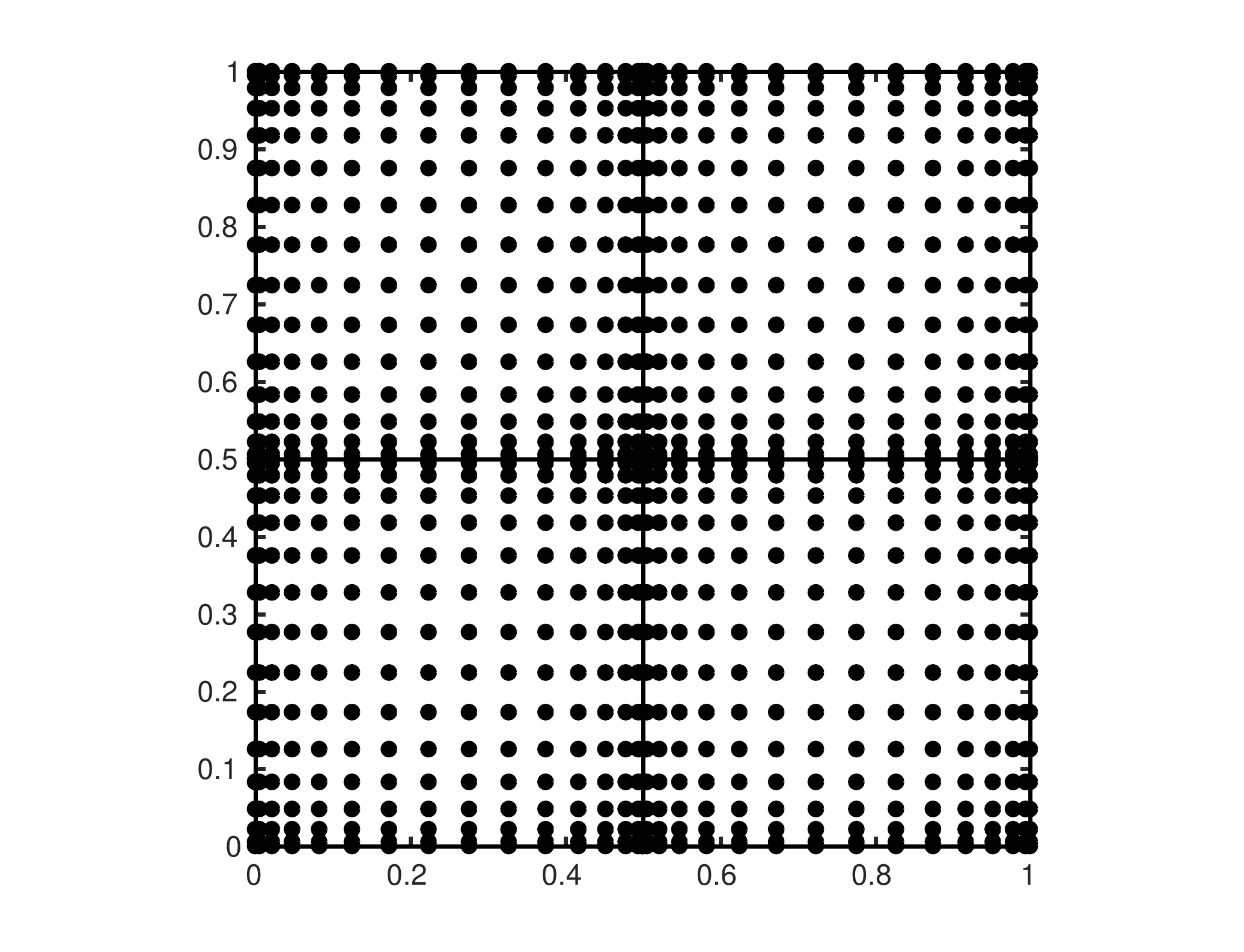}}
\put(81,-3){(b)}
\end{picture}
\caption{\label{fig:interp} Illustration of the 
interpolation points on (a) $\Omega = [0,1]^2$ and (b) its four grandchildren.}
\end{figure}

First, given $n_c$, a tensor product grid of $n_c^2$ Chebychev points is placed
on $\Omega$ and each of its four grandchildren boxes (boxes 4, 5, 6 and 7 in Figure 
\ref{fig:tree_numbering}).  Let $X^{\Omega} = \{\pxx^{\Omega}_l\}_{l=1}^{{n_c}^2}$ 
denote the set of interpolation points defined on box $\Omega$.  
Likewise, let $X^j = \{\pxx^{j}_l\}_{l=1}^{{n_c}^2}$ for $j= 4,5,6,$ and $7$ 
denote the set of interpolation points in box $j$.  Set 
$X^{\rm grand} =\displaystyle \bigcup_{j=4}^7 X^j$. Figure \ref{fig:interp}
illustrates the interpolation points on $\Omega$ and the four
grandchildren when $n_c = 16$ and $\Omega = [0,1]^2$.  

Let $\vct{f}^\Omega$ denote the vector whose entries correspond 
to $f(\pxx)$ evaluated at the points in $X^{\Omega}$ and 
$\mtx{L}_{\rm etk}$ denote the interpolation operator which maps data 
from $X^\Omega$ to $X^{\rm grand}$.  (The notation ${\rm etk}$ stands 
for ``elder to kids.'')   Then $\vct{f}_{\rm app} = \mtx{L}_{\rm etk}\vct{f}^\Omega$ 
is the approximate value of $f(\pxx)$ at the points in $X^{\rm grand}$
interpolated from the values of $f(\pxx)$ at the points in $X^\Omega$.  
Let $\vct{f}_{\rm kids}$ denote the vector whose entries correspond 
to $f(\pxx)$ evaluated at the points in $X^{\rm grand}$.  Let 
$$E_{\rm interp} = \frac{\|\vct{f}_{\rm kids} - \vct{f}_{\rm app}\|_2}{\|\vct{f}_{\rm kids}\|_2}$$
denote the \textit{relative interpolation error}.  If $E_{\rm interp}> \epsilon$, 
$\Omega$ is split into the four grandchildren boxes.  The process is
repeated for each of these smaller boxes.  The process terminates when 
$E_{\rm interp} \leq \epsilon$.

\subsection{Non-uniform HPS solver}
\label{sec:nonuniform}

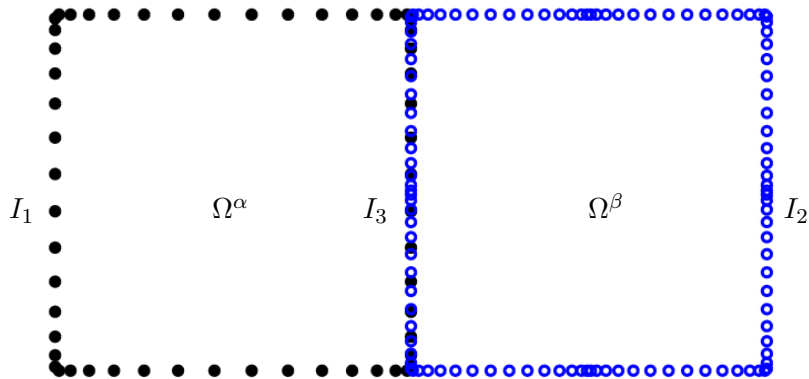
\begin{figure}
\centering
\setlength{\unitlength}{1mm}
\begin{picture}(95,55)
 \put(0,00){\includegraphics[height=55mm]{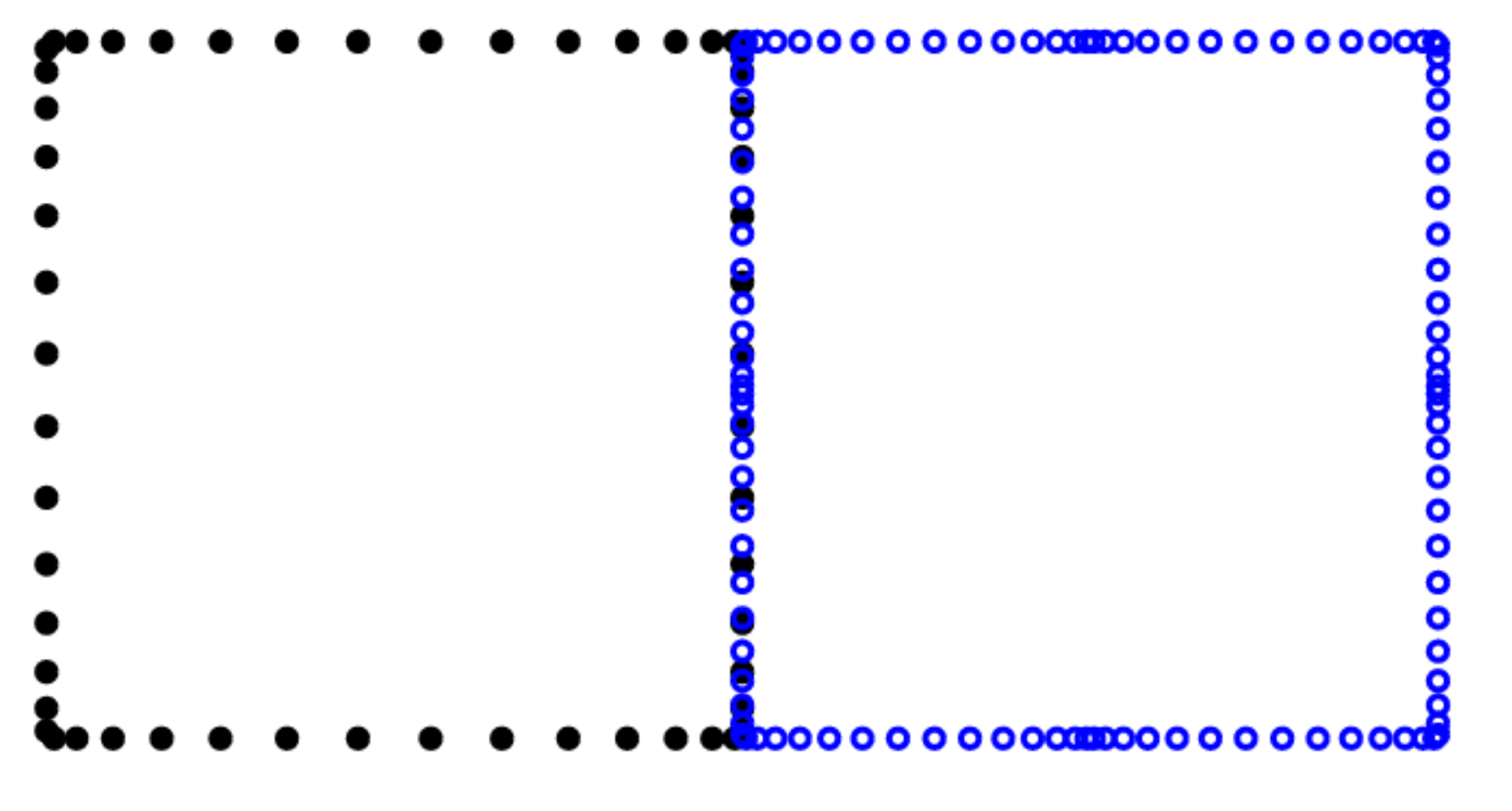}}
\put(24,25){$\Omega^{\alpha}$}
\put(74,25){$\Omega^{\beta}$}
\put(-3,25){$I_{1}$}
\put(100,25){$I_{2}$}
\put(44,25){$I_{3}$}
\end{picture}
\caption{Notation for the merge operation when boxes are on different levels as 
described in Section \ref{sec:nonuniform}.
The rectangular domain $\Omega$ is formed by two squares $\Omega^{\alpha}$ 
and $\Omega^{\beta}$.  The points on the boundary of $\Omega^\alpha$ are 
solid black dots while the points on the boundary of the children of $\Omega^\beta$
are blue hollow dots.}
\label{fig:merge_interp}
\end{figure}

The mesh that results from the adaptive interpolation scheme is likely to be 
highly non-uniform.  While the leaf level operations of the HPS method can
remain the same as for the uniform mesh, the merge operation needs to be 
modified.  Specifically, the boundary operators on the shared interface $I_3$ 
in Figure \ref{fig:merge_interp} need to ``align.''  We chose to align 
the operator via interpolation.  

\begin{figure}[h] 
\centering
\setlength{\unitlength}{1mm}
\begin{picture}(95,55)
\put(-10,00){\includegraphics[height = 50mm]{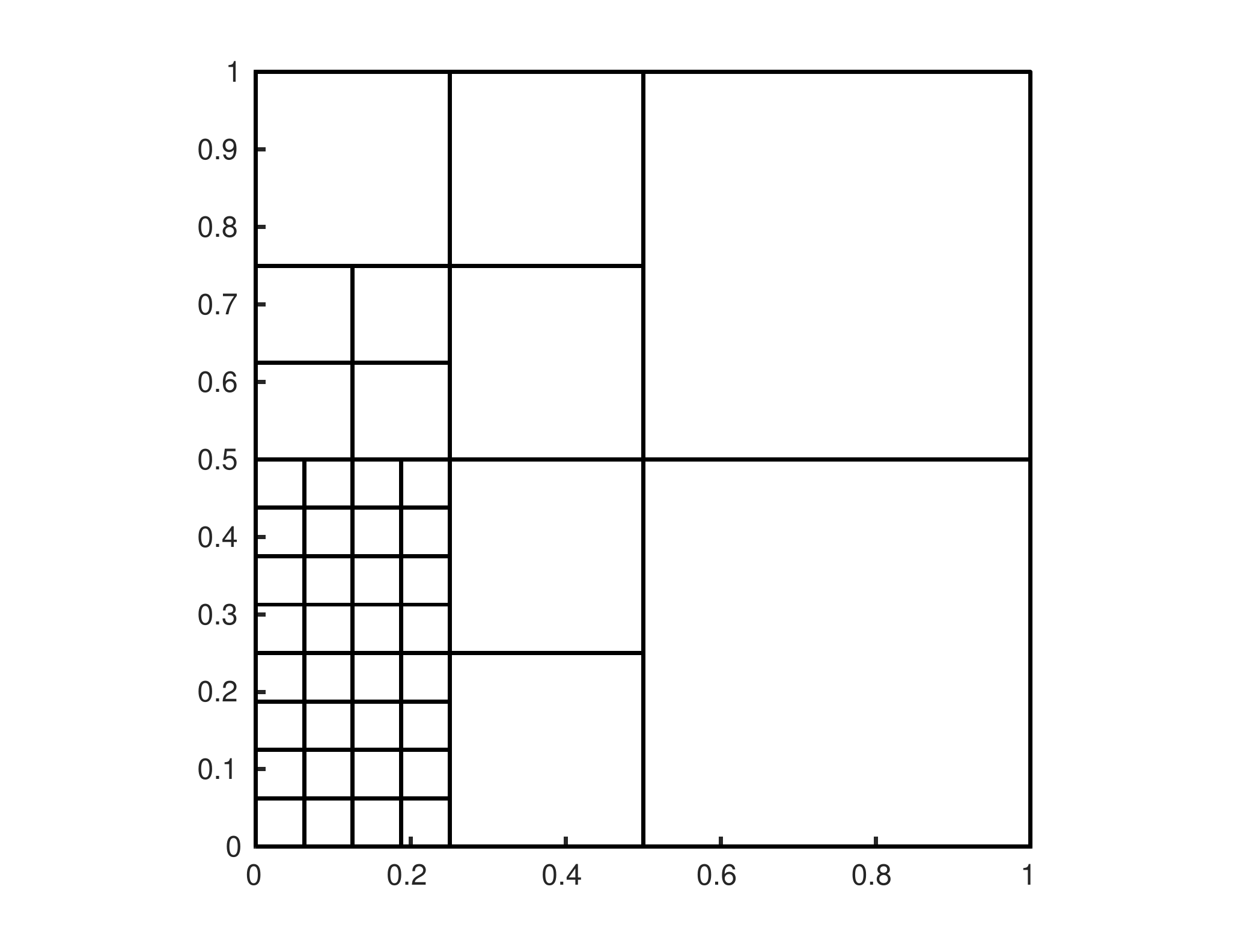}}
\put(21,-3){(a)}
\put(50,00){\includegraphics[height = 50mm]{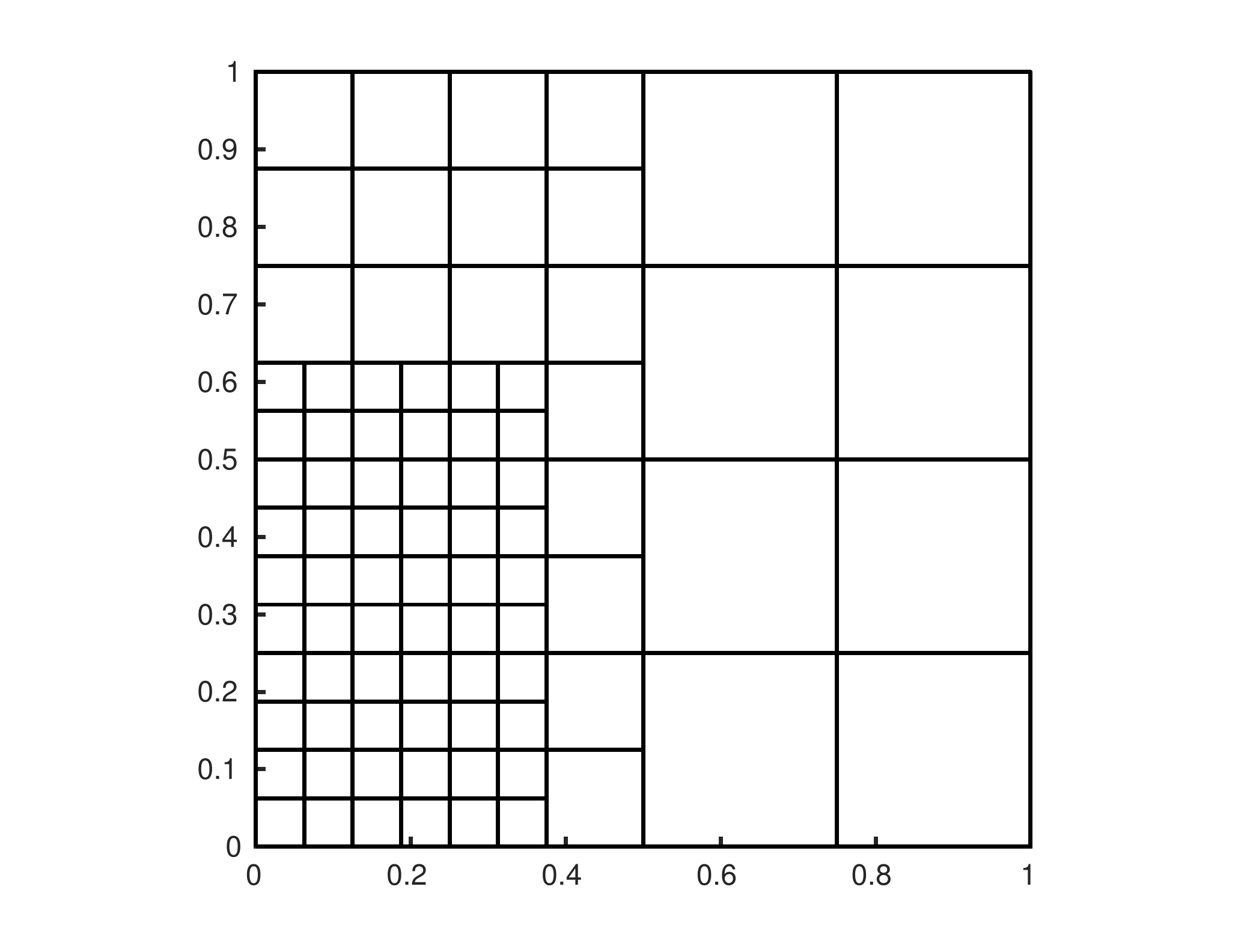}}
\put(81,-3){(b)}
\end{picture}
\caption{\label{fig:levelrestrict} Illustration of the 
mesh resulting from adaptive interpolation applied to 
$f(\pxx) = f(x,y) = e^{-(1000(x-0.11)^2+100(y-0.27)^2)}$ (a) without and (b) with level restriction.
The tolerance was set to $\epsilon = 1e-6$ and $n_c = 16$.}
\end{figure}

Interpolating a super fine mesh to a coarse mesh can be unstable. 
One approach to avoid stability problems is to use nested
interpolation operators that recursively map two panels worth of interpolation points 
to one panels worth of interpolation points.  Alternatively, a level restricted
tree which requires all neighboring boxes be no more than two times bigger
than each other is also stable.  For two dimensional problems, we found the 
constant pre-factors favorable toward the latter approach.  For three dimensional
problems, the nested interpolation will likely be more efficient.
Figure \ref{fig:levelrestrict} illustrates the mesh resulting from the adaptive 
interpolation scheme with and without level restriction.



%
%
%
%
%

%

 
The process of merging two boxes on different levels is straight forward.
For simplicity of presentation, we present the technique for 
merging a leaf box $\alpha$ with a box $\beta$ whose grandchildren 
are leaf boxes.  In this situation, $\Omega^\beta$ has twice as many 
points on its boundary as $\Omega^\alpha$.  Likewise, the DtN matrix 
$\mtx{T}^\beta$ is twice the size of $\mtx{T}^\alpha$.  
Figure \ref{fig:merge_interp} illustrates discretization points on the two boxes.  
The points in $I_3$ from box $\alpha$ do not match the points in 
$I_3$ from box $\beta$.  In order to merge the two 
boxes, we use interpolation.  Let $\mtx{L}_{2t1}$ and $\mtx{L}_{1t2}$ 
denote the interpolation operators that map two panels to one panel on the 
same interval and vice versa.  Since there are $n_c-2$ points on each panel,
the interpolation operators are $n_c-3$ order.  Then the solution and DtN matrices on 
$\Omega^\tau = \Omega^\alpha \cup \Omega^\beta$ are given by 
inserting the interpolation operators into the appropriate locations 
in equations (\ref{eq:Psi}) and (\ref{eq:DtN});
$$\mtx{\Psi}^\tau = \bigl(\TT^{\alpha}_{3,3} - \mtx{L}_{2t1}\TT^{\beta}_{3,3}\mtx{L}_{1t2}\bigr)^{-1}
\bigl[-\TT^{\alpha}_{3,1}\ \big|\ \mtx{L}_{2t1}\TT^{\beta}_{3,2}]$$
and 
$$
\mtx{T}^{\tau} =
\left[\begin{array}{ccc}
\TT_{1,1}^{\alpha} & \mtx{0} \\
\mtx{0} & \TT_{2,2}^{\beta }
\end{array}\right] +
\left[\begin{array}{c}
\TT_{1,3}^{\alpha} \\
\TT_{2,3}^{\beta}\mtx{L}_{1t2}
\end{array}\right]\,
\mtx{\Psi}^\tau.
$$

\begin{remark}
 Merging two boxes on different levels was also presented in 
 \cite{2016_bodyload}.  There the DtN operator $\TT^\beta$ is 
 interpolated to a corresponding operator so that the number
 of points per edge matches the operator on box $\alpha$.  
 The method presented in this paper does not 
 take this approach since the coarse sampling of boundary data on $\Omega^1=\Omega$ 
 may not be sufficient resulting in a loss of accuracy.
\end{remark}

\subsection{Indicator for refinement}
\label{sec:mark}
This section presents a technique for identifying which leaf boxes
need further refinement.  The technique utilizes the fact that the
local problem can be fully resolved even though the error is large
due to incorrect boundary data.  To determine if the local problem
is fully resolved, the method we propose looks at the tail coefficients
of the approximate solution written in a Chebychev expansion on the 
leaf boxes.  This technique is inspired by \cite{lee_greengard_stiff}
where a decay condition on the local Chebychev expansion coefficients was used to build an 
adaptive discretization technique for one dimensional integral equations.  

Since each leaf box has a modified tensor product basis, the approximate 
solution $u_{\rm app}$ at any discretization point can be written as the product of two Lagrange
polynomials $\Phi(x)$ and $\Psi(y)$ with $x$ and $y$ interpolation nodes respectively.  
Thus we will look at the 
\textit{directional} Chebychev coefficients to build a refinement criterion.
Recall that for one dimensional interpolation 
the Lagrange interpolant through Chebychev points can 
be expressed as a partial Chebychev expansion with coefficients
that can be found via the Fast Fourier Transform (FFT)\cite{2000_boyd}.  
If the basis is sufficient to capture the solution locally, the series
will be convergent and thus we can approximate the contributions from the
remainder of the series by looking at the last few coefficients of the 
local expansions \cite{lee_greengard_stiff}.

Specifically, for a leaf box $\Omega^\tau$, 
let $\{x_j\}_{j=1}^{n_c}$ and $\{y_j\}_{j=1}^{n_c}$ 
denote the one dimensional Chebychev interpolation points that the 
$x$ and $y$ coordinates of any discretization point.  For a fixed $x_j$, 
$j =2,\ldots, n_c-1$, let the vector $\vct{B}_{j}$ of length $n_c$ denote the Chebychev 
coefficients in the $y-$direction of the approximate solution along the 
line $x = x_j$ for $y \in[y_1,y_n]$; i.e. the entries in $\vct{B}_j$ correspond to the
coefficients of $\Psi(y)$ written in terms of Chebychev polynomials along the line $x=x_j$ in $\Omega^\tau$.  
We define  
$$S_y = \max_{j=2,\ldots,n_c-1} \left(|\vct{B}_j(n_c-1)|+|\vct{B}_j(n_c)-\vct{B}_j(n_c-2)|\right)$$
be the indicator of the decay of the coefficients in the $y-$direction.  
Likewise for a fixed $y_j$, $j=2,\ldots,n_c-1$, 
let $\vct{C}_j$ denote vector of length $n_c$ containing the Chebychev coefficients of 
the approximate solution along the line $y=y_j$ in $\Omega_\tau$ (i.e. the $x-$direction 
coefficients) and define
$$S_x = \max_{j=2,\ldots,n_c-1} \left(|\vct{C}_j(n_c-1)|+|\vct{C}_j(n_c)-\vct{C}_j(n_c-2)|\right)$$ 
to be the indicator of the decay of the coefficients in the $x-$direction.  
Then, for each leaf box $\tau$, we define
$S_\tau = \max\{S_x,S_y\}$.  This yields a measure for how well the local basis
is able capture solutions to the partial differential equation restricted to $\Omega^\tau$. 
Let 
$$ S_{\rm div} = \frac{1}{4}\displaystyle\max_{{\rm leaf \ boxes} \ \tau} S_\tau.$$
This gives a measure for how accurate we should hope the tail of the local expansions should 
be. Any leaf box that does not meet this requirement, i.e. $S_\tau >S_{\rm div}$, 
is marked for further refinement.  Algorithm 3 presents a 
pseudocode for determining which leaf boxes need refinement.


\begin{figure}[h!]
\fbox{
\begin{minipage}{120mm}
\begin{center}
\textsc{Algorithm 3} (Refinement indicator)
\end{center}

This algorithm presents the technique for determining which leaf boxes need 
additional refinement. It assumes a tree structured mesh and the corresponding 
direct solver are given.

\rule{\textwidth}{0.5pt}

\begin{tabbing}
\mbox{}\hspace{7mm} \= \mbox{}\hspace{6mm} \= \mbox{}\hspace{6mm} \= \mbox{}\hspace{6mm} \= \mbox{}\hspace{6mm} \= \kill
(1)\> \textbf{for} $\tau = N_{\rm boxes},\,N_{\rm boxes}-1,\,N_{\rm boxes}-2,\,\dots,\,1$\\
(2)\> \> \textbf{if} ($\tau$ is a leaf)\\
(3)\> \> \> \textbf{for} $j = 2,\ldots, n_c-1$\\
(4)\> \> \> \> Compute the $y$-directional Chebychev coefficients $\vct{B}_j$ of the \\
\> \> \> \>approximate solution on $\tau$.\\
(5) \> \> \> \textbf{end for}\\
(6) \> \> \> Let $\displaystyle S_y = \max_{j=2,\ldots,n_c-1} \left(|\vct{B}_j(n_c-1)|+|\vct{B}_j(n_c)-\vct{B}_j(n_c-2)|\right)$\\
(7)\> \> \> \textbf{for} $j = 2,\ldots, n_c-1$\\
(8)\> \> \> \> Compute the $x$-directional Chebychev coefficients $\vct{C}_j$ of the \\
\> \> \> \> approximate solution on $\tau$.\\
(9) \> \> \> \textbf{end for}\\
(10) \> \> \> Let $\displaystyle S_x = \max_{j=2,\ldots,n_c-1} \left(|\vct{C}_j(n_c-1)|+|\vct{C}_j(n_c)-\vct{C}_j(n_c-2)|\right)$\\
(11) \> \> \> Let $S_\tau = \max\{S_x,S_y\}$\\
(12) \> \> \textbf{end if}\\
(13) \> \textbf{end for}\\
 (14) \> Let $ S_{\rm div} = \frac{1}{4}\left(\displaystyle\max_{{\rm leaf \ boxes }\ \tau} S_\tau\right).$\\
(15)\> \textbf{for} $\tau = N_{\rm boxes},\,N_{\rm boxes}-1,\,N_{\rm boxes}-2,\,\dots,\,1$\\
(16)\> \> \textbf{if} ($\tau$ is a leaf)\\
(17)\> \>  \> \textbf{if} $S_\tau > S_{\rm div}$\\
(18) \> \> \> \> Add $\tau$ to the refinement list.\\
(19) \> \>  \> \textbf{end if} \\
(20)\> \textbf{end for}
\end{tabbing}
\end{minipage}}
\label{alg:adapt_disc}
\end{figure}

\subsection{Updating the solver}
\label{sec:update}
Once the list of leaf boxes marked for refinement is made, we need to solve
(\ref{eq:basic}) with the refined grid to determine if the mesh gives the 
desired accuracy. Constructing the direct 
solver from scratch is computationally expensive and unnecessary when the 
refinement is localized.  This section presents a technique for building
the solver for the refined mesh while making use of the existing solver.  The 
key observation is the fact that the solution
technique is naturally domain decomposing.  This means that the only parts 
of the solver that need to be modified are the parts that touch the refined regions.

The first step in this process is to make a list of all boxes affected by the 
local refinement.  To do this, starting from the list of boxes 
refined, we sweep the binary tree making note of all the 
ancestors affected.  For example, if boxes 16, and 18 were the only boxes
marked for refinement in Figure \ref{fig:tree_numbering}, 
the solver would need to update the operators for boxes 1, 2, 4, 8 and 9.  
The operators for the other boxes need not be touched.

Next DtN and solution matrices are constructed by moving through 
the list of effected boxes starting from the bottom of the tree (i.e.
first processing the leaf boxes then its ancestors in order of 
ancestry).

\begin{remark}
 Further acceleration can be gained by creating new tree structures based 
 on the refinement regions.  For the example where boxes 16 and 18 in Figure
 \ref{fig:tree_numbering} require refinement, DtN and solution matrices can 
 be constructed for the union of boxes 5, 6, 7, 17 and 19.  Then 
 the computation is limited to the boxes 16, 18, their union and gluing the
 union with the remainder of the geometry.  For the problems under consideration
 in this manuscript, this technique was not employed. 
\end{remark}

\section{Numerical results}
\label{sec:numerics}
This section illustrates the ability of the adaptive discretization technique
to solve a collection of problems.
First, in section \ref{sec:elliptic}, 
three problems suggested in \cite{MITCHELL2013} to test adaptive discretization 
techniques for elliptic PDEs are considered.  For each of these problems
the solution is known but each poses a different challenge for adaptive 
discretization techniques.   
Section \ref{sec:helm} considers two Helmholtz problems: a low- to mid-frequency
constant coefficient problem with a source and a high frequency variable 
coefficient problem.  The globally oscillatory
nature of the solution adds to the challenge of accurately discretizing these
problems.

The following quantities are reported.\\
\begin{tabular}{rl}
	$n_c$: &the choice of discretization order \\
 $N_i$: & the number of leaf boxes after adaptive interpolation \\
 $N_f$: & the number of leaf boxes after adaptive discretization\\
 $T_i$: & the time in seconds for the adaptive interpolation step\\
 $T_f$: & the time in seconds for the adaptive discretization step\\
 $T_s$: & the time in seconds to apply the resulting solver \\
$R$: & the memory in GB for storing the direct solver \\
\end{tabular}

For all experiments, the uniform discretization technique is applied
for comparison purposes.  $T_{\rm pre}$ is used to report the 
time in seconds for discretizing the PDE and building the direct solver.


To report on the accuracy of the solution techniques, we report 
$$E_{\rm rel} = \frac{1}{N_f} \sum_{{\rm leaf \ boxes} \tau} E^\tau_{\rm rel}$$
where $$E^\tau_{\rm rel} = \frac{\|\vct{u}^\tau_{\rm app}-\vct{u}^\tau_{\rm ref}\|_2}{\|\vct{u}^\tau_{\rm ref}\|_2}$$
for leaf box $\tau$, $\vct{u}^\tau_{\rm app}$ is the approximate solution at the discretization points
on $\tau$, and $\vct{u}^\tau_{\rm ref}$ is the reference solution evaluated at the discretization points
on $\tau$.  For the problems where the solution is known, the reference solution is the exact solution.
For problems where the solution is unknown, the reference solution is given by an approximate solution 
obtained by running the uniform HPS method until convergence.


\subsection{Problems with known solutions}
\label{sec:elliptic}
This section reports the performance of the solution techniques for three 
problems where the solution is known and the partial differential equation
has smooth coefficients on the domain $\Omega = (0,1)^2$.  The problems 
under consideration are the following:\\
\noindent
\textit{Boundary layer:}  The Dirichlet boundary value problem 
\begin{equation*}
    -\alpha\nabla^2u + 2\frac{\partial u}{\partial x} + \frac{\partial u}{\partial y} = f(x,y),
\end{equation*}
where the solution is given by 
\begin{equation*}
    u(x,y) = (1-e^{-(1-x)/\alpha})(1-e^{-(1-y)/\alpha})\cos(\pi(x+y))
\end{equation*}
and the parameter $\alpha= 10^{-3}$ determines the steepness of the boundary layer. \\
\noindent 
\textit{Locally oscillatory solution:} The Dirichlet boundary value problem 
\begin{equation*}
    -\nabla^2u - \frac{1}{(\alpha+\sqrt{x^2+y^2})^4}u = f(x,y),
\end{equation*}
 where the solution is given by 
\begin{equation*}
    u(x,y) = \sin\left(\frac{1}{\alpha+\sqrt{x^2+y^2}}\right)
\end{equation*}
and the parameter $\alpha = \frac{1}{10\pi}$ determines the number of oscillations
in the solution.  The oscillations are clustered near the origin. \\
\noindent
\textit{Wave front:} The Poisson Dirichlet boundary value problem
where the solution is given by 
\begin{equation*}
    u(x,y) = \tan^{-1}(50(\sqrt{(x+0.05)^2 + (y+0.05)^2}-0.7)).
\end{equation*}

\begin{figure}[h!]
\centering
 \begin{tabular}{ccc}
  \includegraphics[height = 45mm]{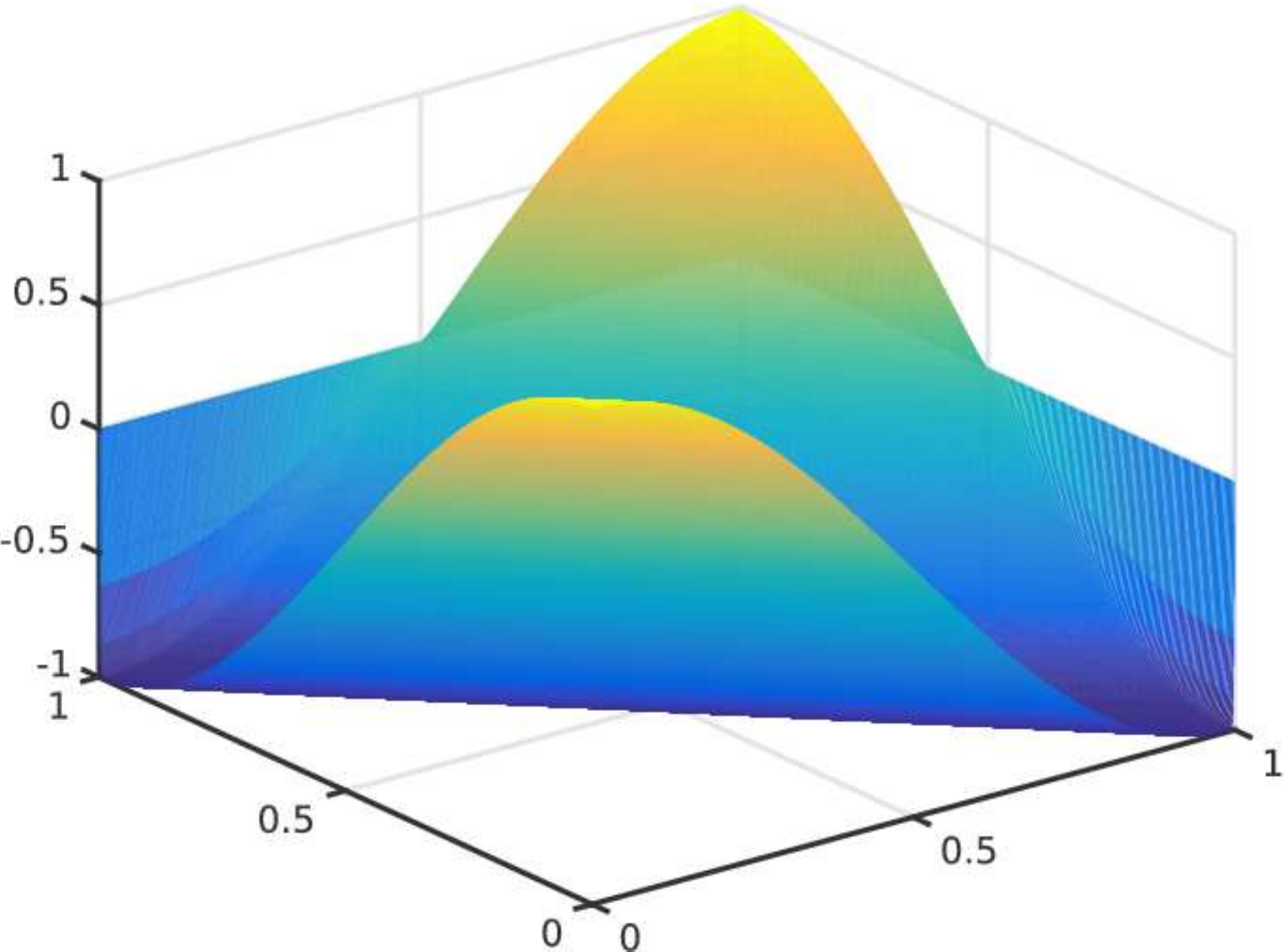} & \mbox{\hspace{20mm}} & 
    \includegraphics[height = 45mm]{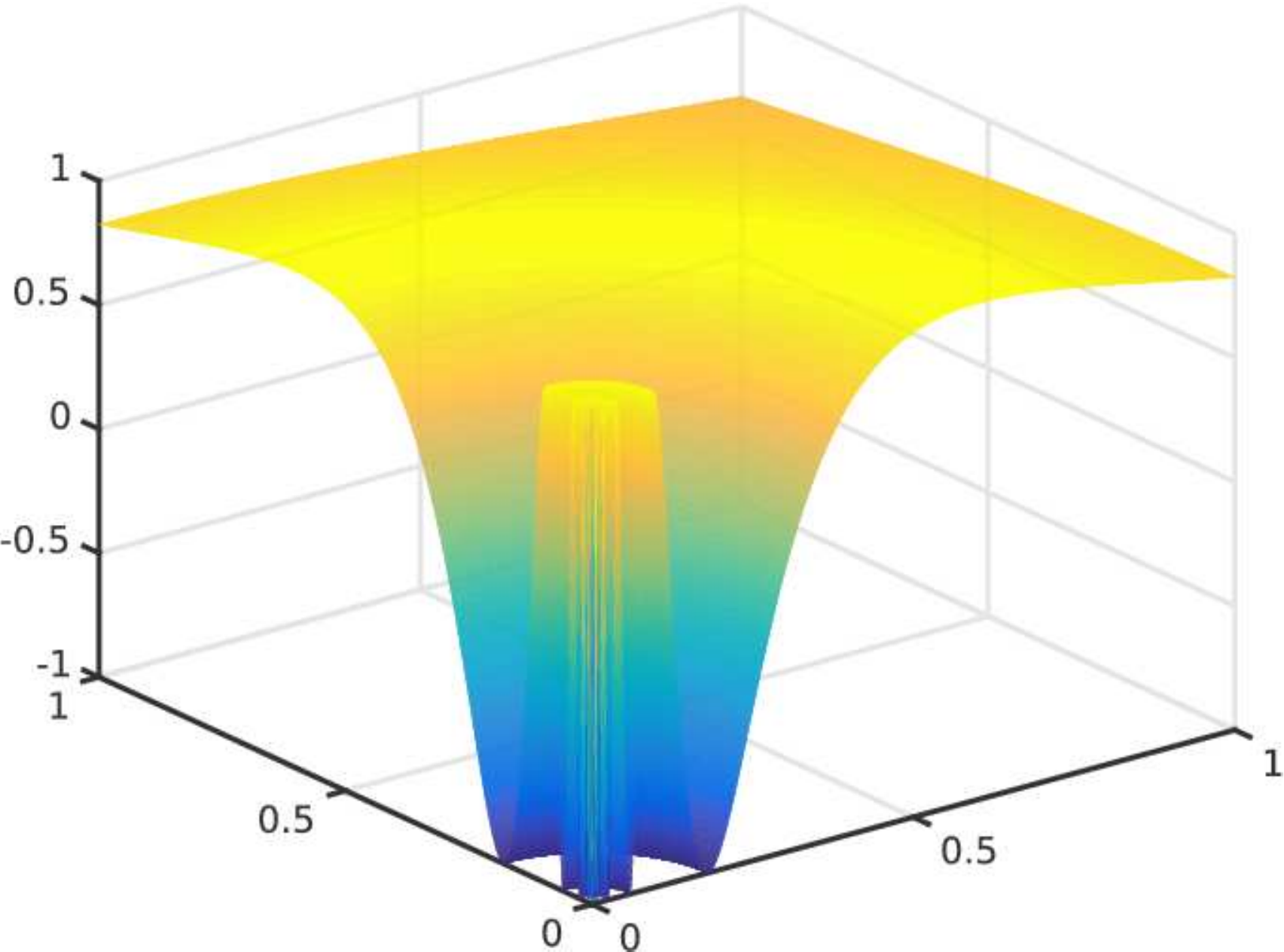}\\
    &&\\
    (a) Boundary layer & & (b) Locally oscillatory \\
    \multicolumn{3}{c}{\includegraphics[height = 45mm]{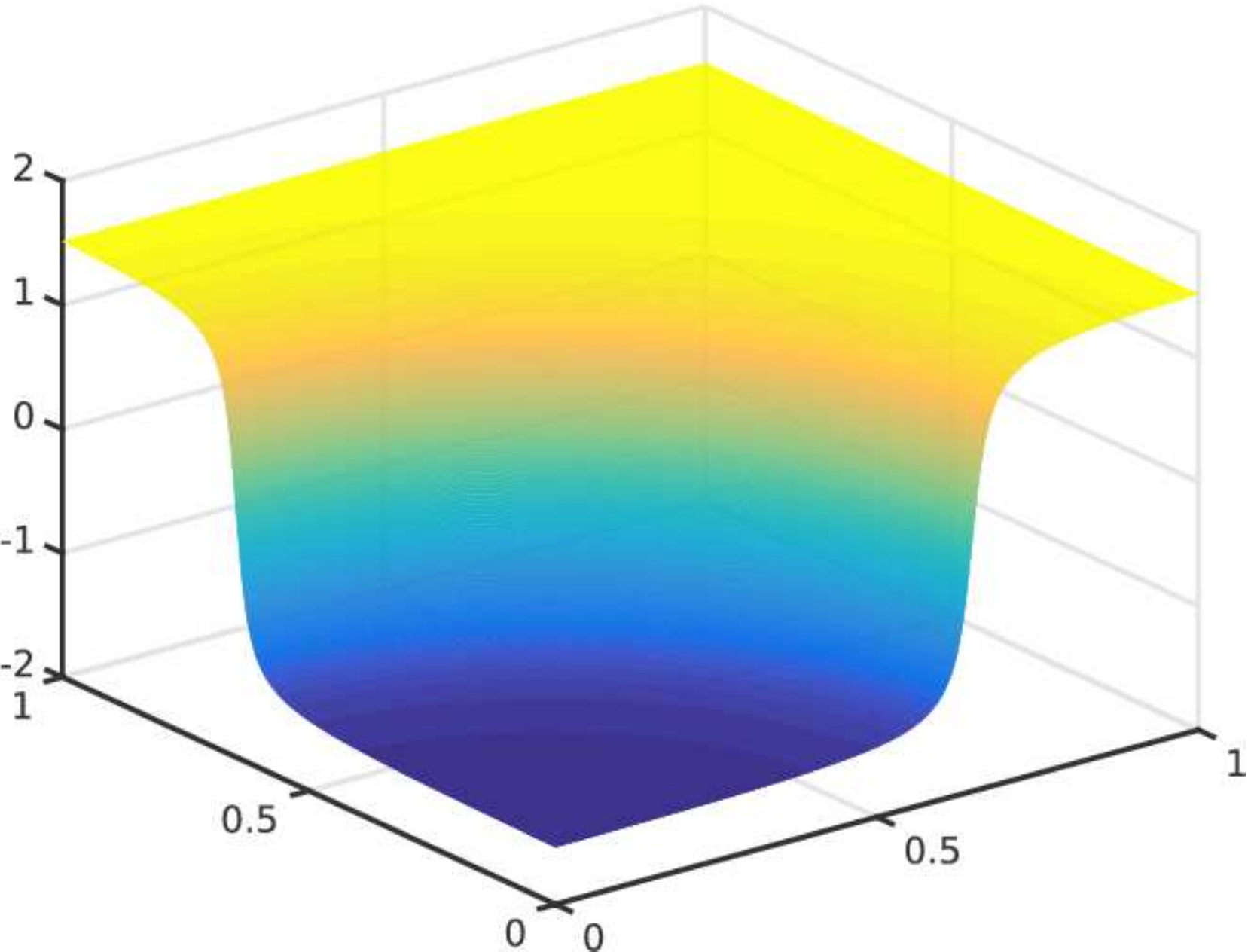}}\\
       &&\\
    \multicolumn{3}{c}{(c) Wave front}
 \end{tabular}
\caption{\label{fig:solns}  Illustration of the solutions to the problems under 
consideration in section \ref{sec:elliptic}: (a) Boundary layer problem, (b)
Problem with a locally oscillatory solution, and (c) Problem where the solution is
a wavefront.}
\end{figure}

Figure \ref{fig:solns} illustrates the solutions to each of these problems.  
Table \ref{tab:DtN} reports the performance of the method for each of these
problems with the stopping tolerance set to $\epsilon  = 10^{-5}$ and different
discretization orders $n_c = 8, 16$ and $32$.   For all of the experiments the 
adaptive algorithm achieves the desired tolerance.  In fact, for most 
of the experiments the discretization technique achieves better than the
desired tolerance.  
The results also indicate that since the solutions to the 
boundary layer and locally oscillatory problem are ``nicer''
than the coefficients of the partial differential equation, 
the mesh achieved via the adaptive interpolation 
technique is more than sufficient for resolving the problem.  
For the wave front problem discretized
with the low order basis, the adaptive discretization technique is 
needed to achieve the user specified tolerance.  For all the experiments, it 
is computationally beneficial (less expensive to achieve the same or 
better accuracy) to use a higher order discretization.  The timing 
results re-enforce the benefit of using the high order basis.

\begin{figure}[h!]
\centering
 \begin{tabular}{ccc}
  \includegraphics[height = 45mm]{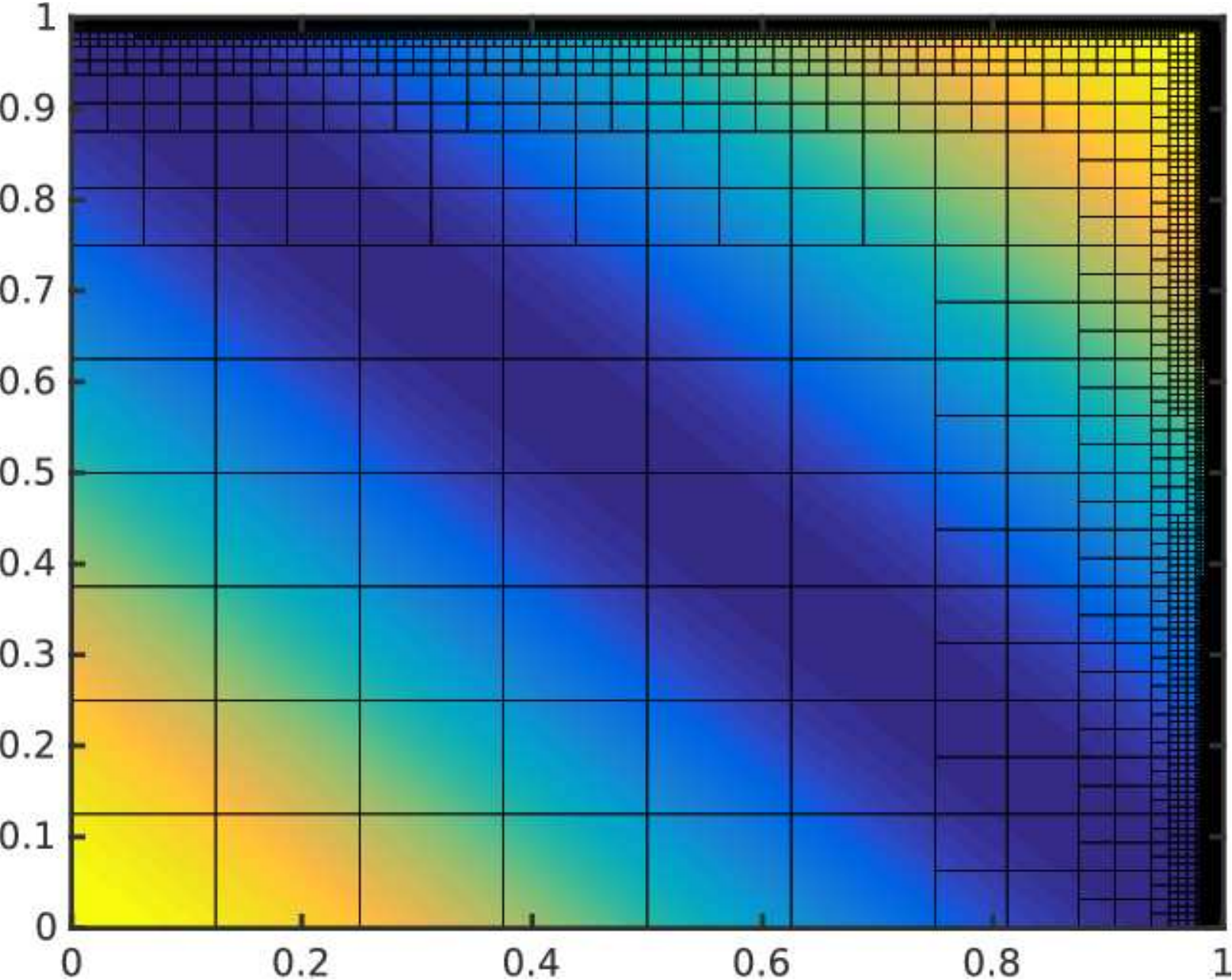} & \mbox{\hspace{15mm}} & 
    \includegraphics[height = 45mm]{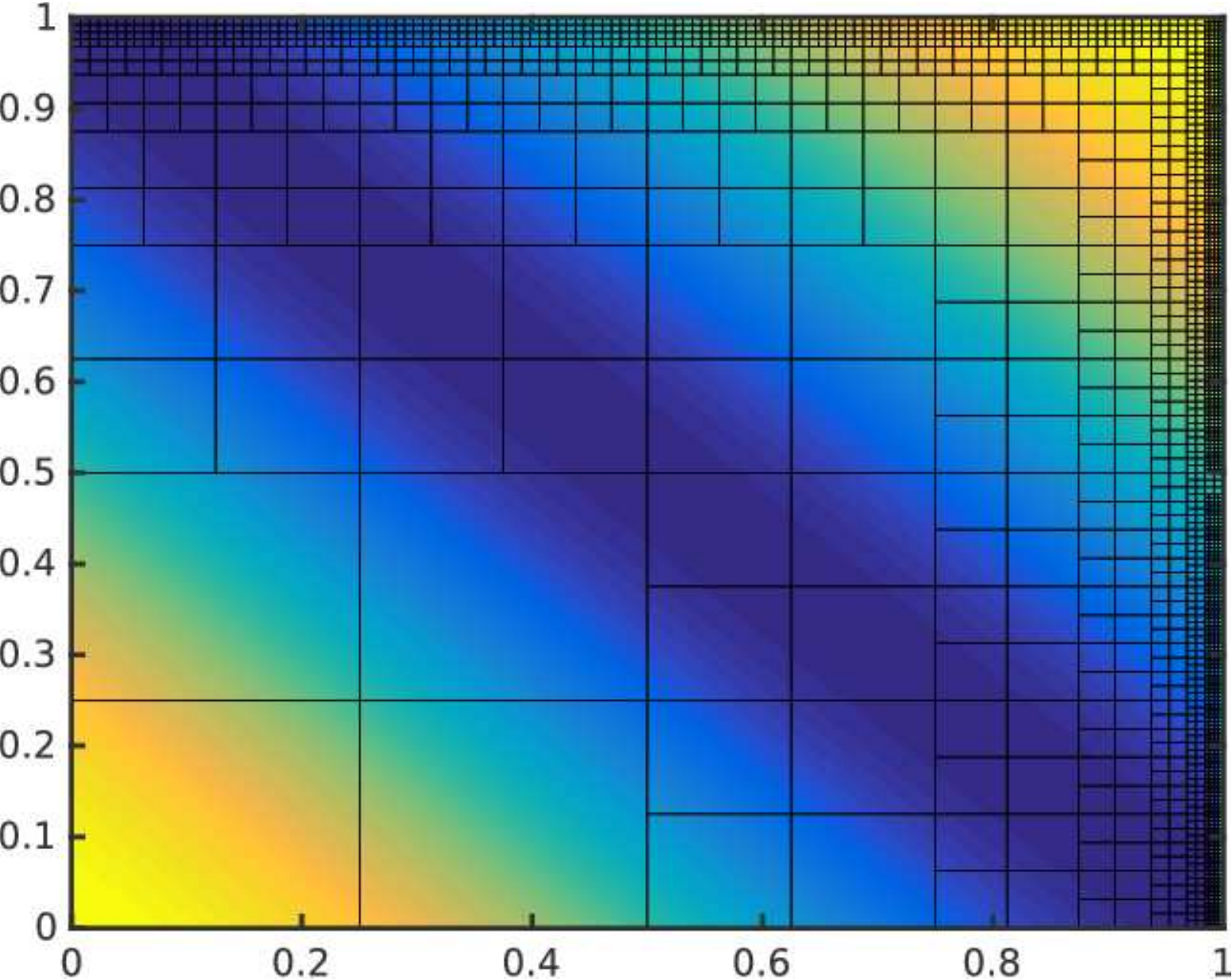}\\
    (a) Boundary layer with $n_c= 8$ & & (b) Boundary layer with $n_c= 16$\\
  \includegraphics[height = 45mm]{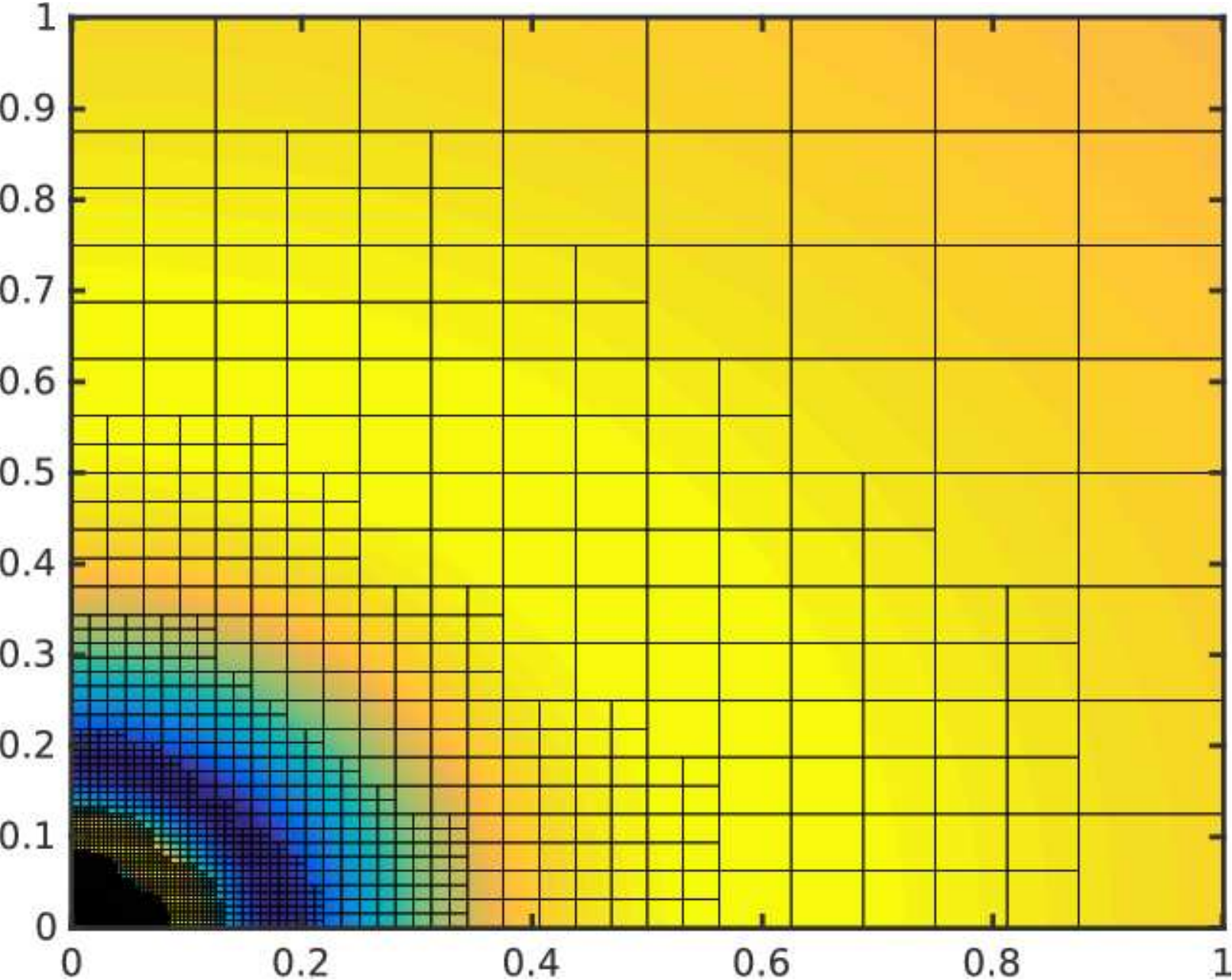} & \mbox{\hspace{15mm}} & 
    \includegraphics[height = 45mm]{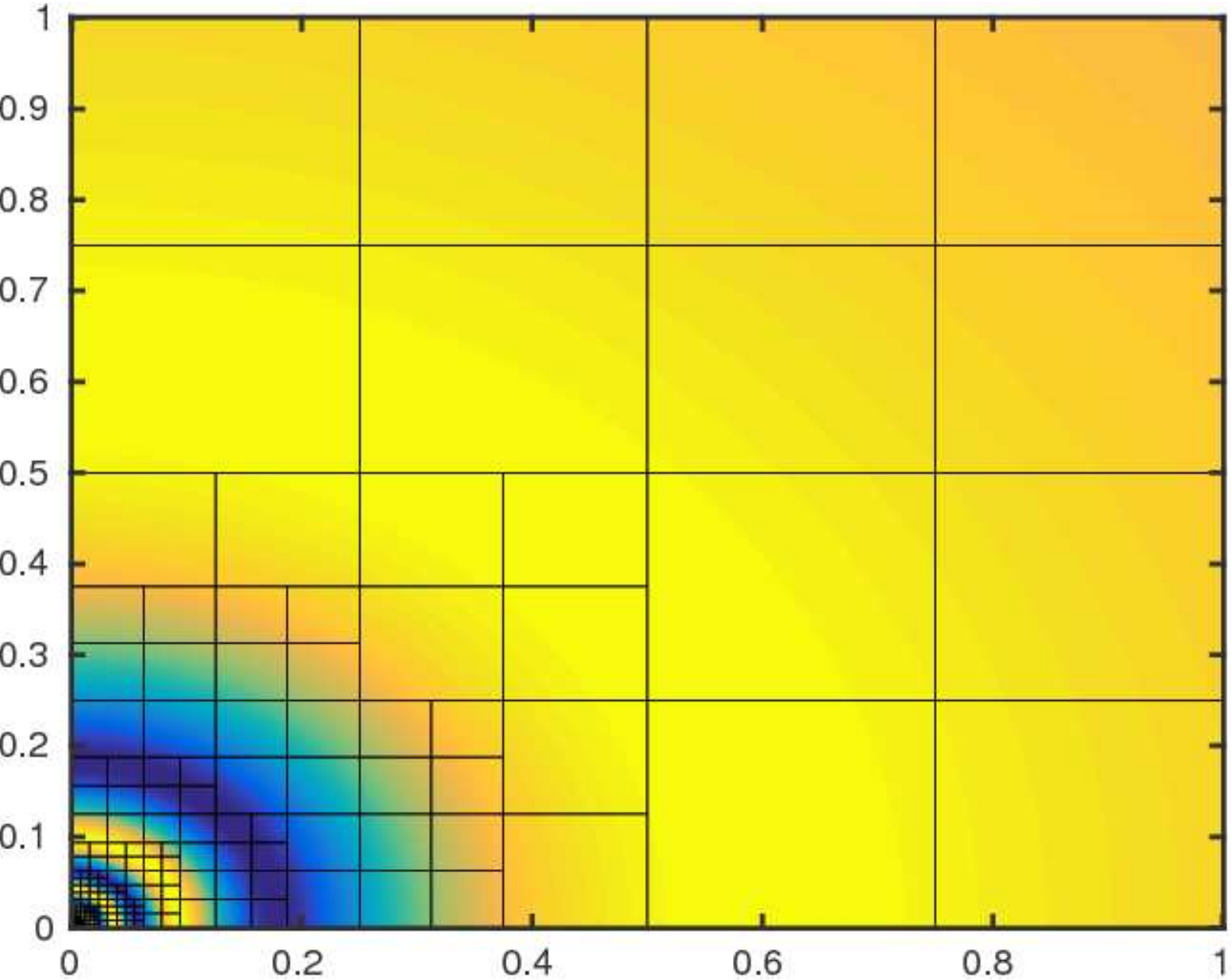}\\
    (c) Oscillatory solution with $n_c= 8$ & & (d) Oscillatory solution with $n_c= 16$ \\
  \includegraphics[height = 45mm]{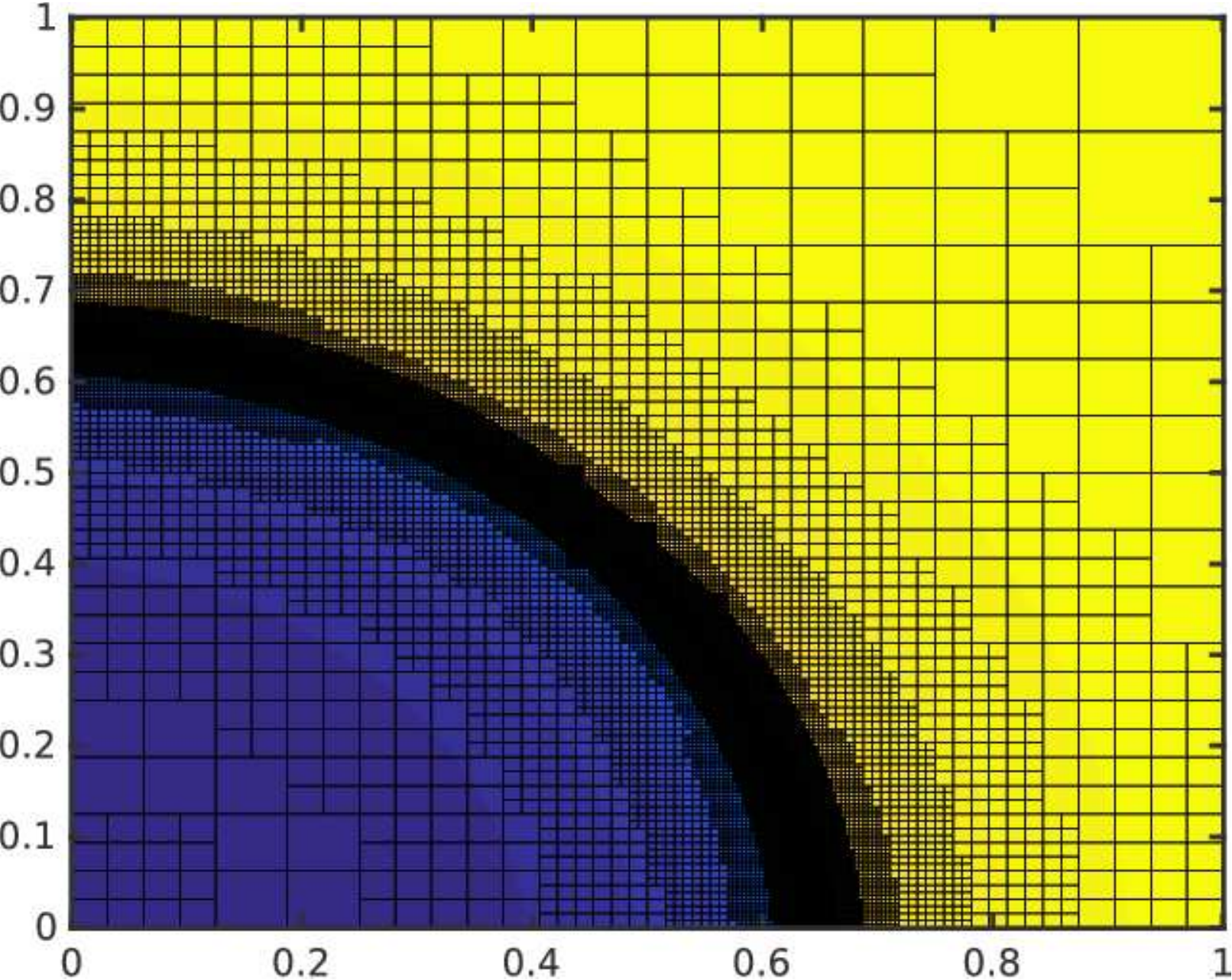} & \mbox{\hspace{15mm}} & 
    \includegraphics[height = 45mm]{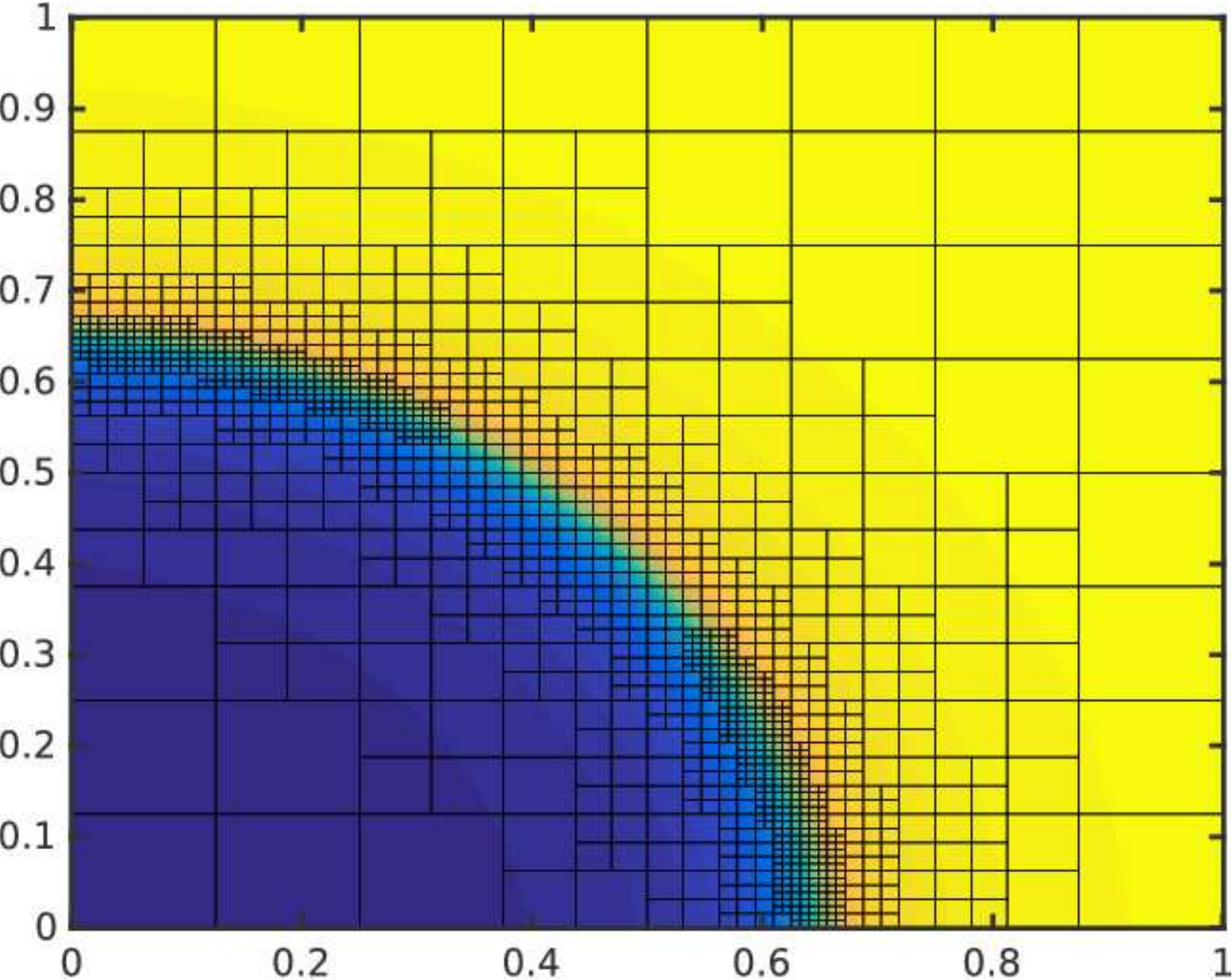}\\
    (e) Wave front solution with $n_c= 8$ & & (f) Wave front solution with $n_c= 16$ \\
 \end{tabular}
\caption{\label{fig:mesh1}  Illustration of the mesh generated by the 
adaptive discretization technique (overlayed on the solution) with $n_c = 8$ and $16$
for the problems in section \ref{sec:elliptic}.}
\end{figure}

Figure \ref{fig:mesh1} illustrates the mesh overlayed on the solution for each 
experiment.  The mesh shows that the method is finding the areas where refinement
is necessary.  The denseness of the leaf boxes visualize the additional cost of 
using a low order method.

Tables \ref{tab:DtN_layer}-\ref{tab:DtN_wave} report the time in seconds for 
applying the uniform discretization technique to the three partial differential 
equations $T_f$ and applying the direct solver $T_s$
for different orders of discretization $n_c$ and numbers of leaf boxes $N_f$.  
The memory $R$ for storing the direct solver and the relative error $E_{\rm rel}$ are 
also reported.  The cost for building the direct solver is more expensive for the 
adaptive discretization technique than for the uniform discretization.  This is 
the case for most adaptive discretization techniques.  The cost of applying the 
direct solver from the adaptive discretization is less expensive than applying the
solver from the uniform discretization.  To achieve the same accuracy as the 
adaptive discretization, the uniform discretization requires more leaf boxes 
and is more expensive (measured by adding the cost testing the finer grids).

\begin{table}[h!]
\centering
\begin{tabular}{|c|c|c|c|c|c|c|c|c|}
	\hline
	Problem & $n_c$ & $N_i$ & $N_f$ & $T_i$ & $T_f$ & $T_s$ & $R$ & $E_{rel}$ \\
	\hline
	\multirow{3}{*}{Boundary layer}      
	& 8  & 66610 &  66610 & 3.58e$+$03 & 1.99e$+$02 & 1.01e$+$01 & 33.4 & 5.39e$-$09 \\ \cline{2-9}                   
	& 16 &  2194 &   2194 & 1.60e$+$01 & 3.18e$+$01 & 1.29e$+$00 & 2.99 & 7.27e$-$10 \\ \cline{2-9}
	& 32 &   316 &    316 & 2.42e$+$01 & 5.70e$+$01 & 4.97e$-$01 & 3.46 & 3.52e$-$13 \\ 
    \hline
	\multirow{3}{*}{Locally oscillatory} 
	& 8  & 21247 &  21247 & 5.08e$+$02 & 4.24e$+$01 & 3.18e$+$00 & 5.57 & 1.35e$-$08 \\ \cline{2-9}
	& 16 &   487 &    487 & 4.43e$+$00 & 5.97e$+$00 & 2.21e$-$01 & 0.78 & 1.93e$-$08 \\ \cline{2-9}
	& 32 &   232 &    232 & 7.05e$+$00 & 2.14e$+$01 & 3.06e$-$01 & 3.15 & 4.06e$-$09 \\ 
    \hline
    \multirow{3}{*}{Wave front}          
    & 8  & 44392 & 148087 & 1.28e$+$03 & 7.24e$+$03 & 2.34e$+$01 & 56.0 & 5.43e$-$04 \\ \cline{2-9}                   
	& 16 &  1405 &   1405 & 1.42e$+$01 & 1.56e$+$01 & 4.96e$-$01 & 1.39 & 4.36e$-$11 \\ \cline{2-9}
	& 32 &   349 &    349 & 1.60e$+$01 & 1.40e$+$02 & 4.60e$-$01 & 3.74 & 5.20e$-$12 \\ 
    \hline
\end{tabular}
 \caption{\label{tab:DtN}  
 Timing, memory and error results for applying the adaptive technique to each of 
 the experiments in section  \ref{sec:elliptic} with 
 different orders of discretization $n_c$.}
\end{table}

\begin{table}[h!]
\centering
\begin{tabular}{|c|r|c|c|c|c|}
	\hline
	$n_c$ & $N_f$ & $T_{\rm pre}$ & $T_s$ & $R$ & $E_{rel}$ \\
	\hline     
	\multirow{7}{*}{16}
	&     4 & 9.84e$-$02 & 4.07e$-$02 & 0.003 & 9.46e$-$01 \\ \cline{2-6}                   
	&    16 & 2.14e$-$01 & 2.01e$-$02 & 0.01  & 2.63e$-$01 \\ \cline{2-6}
	&    64 & 4.39e$-$01 & 2.95e$-$02 & 0.06  & 4.21e$-$02 \\ \cline{2-6}
	&   256 & 1.40e$+$00 & 1.24e$-$01 & 0.24  & 2.19e$-$03 \\ \cline{2-6}
	&  1024 & 5.83e$+$00 & 5.45e$-$01 & 1.05  & 2.79e$-$05 \\ \cline{2-6}
	&  4096 & 2.49e$+$01 & 3.93e$+$00 & 4.52  & 5.40e$-$08 \\ \cline{2-6}
	& 16384 & 1.19e$+$02 & 1.76e$+$01 & 19.3  & 1.93e$-$11 \\
    \hline
    %
    \multirow{5}{*}{32}
	&     4 & 3.69e$-$01 & 3.93e$-$02 & 0.04 & 7.92e$-$02 \\ \cline{2-6}                   
	&    16 & 1.21e$+$00 & 3.61e$-$02 & 0.16 & 5.15e$-$03 \\ \cline{2-6}
	&    64 & 3.89e$+$00 & 1.32e$-$01 & 0.67 & 4.72e$-$05 \\ \cline{2-6}
	&   256 & 1.51e$+$01 & 1.60e$+$00 & 2.76 & 1.74e$-$08 \\ \cline{2-6}
	&  1024 & 1.01e$+$02 & 6.36e$+$00 & 11.4 & 1.22e$-$13 \\
	\hline
\end{tabular}
\caption{\label{tab:DtN_layer}
 Timing, memory and error results for applying the uniform discetization technique to 
 the boundary layer problem with different orders of discretization $n_c$.}
\end{table}
\begin{table}[h!]
\centering
\begin{tabular}{|c|r|c|c|c|c|}
	\hline
	$n_c$ & $N_f$ & $T_{\rm pre}$ & $T_s$ & $R$ & $E_{rel}$ \\
	\hline     
	\multirow{8}{*}{16}
	&      4 & 1.25e$-$01 & 2.79e$-$02 & 0.003 & 4.78e$-$01 \\ \cline{2-6}
	&     16 & 2.20e$-$01 & 1.82e$-$02 & 0.01  & 6.44e$-$01 \\ \cline{2-6}
	&     64 & 5.12e$-$01 & 3.39e$-$02 & 0.06  & 5.97e$-$01 \\ \cline{2-6}
	&    256 & 1.63e$+$00 & 9.50e$-$02 & 0.24  & 1.11e$-$01 \\ \cline{2-6}
	&   1024 & 6.03e$+$00 & 3.78e$-$01 & 1.05  & 6.16e$-$04 \\ \cline{2-6}
	&   4096 & 2.32e$+$01 & 2.52e$+$00 & 4.52  & 2.55e$-$05 \\ \cline{2-6}
	&  16384 & 1.00e$+$02 & 1.50e$+$01 & 19.3  & 3.25e$-$06 \\ \cline{2-6}
	&  65536 & 4.27e$+$02 & 6.51e$+$01 & 82.4  & 4.09e$-$07 \\
    \hline
	\multirow{7}{*}{32}
	&      4 & 4.05e$-$01 & 3.57e$-$02 & 0.04  & 4.88e$-$01 \\ \cline{2-6}
	&     16 & 1.07e$+$00 & 3.61e$-$02 & 0.16  & 6.93e$-$03 \\ \cline{2-6}
	&     64 & 3.27e$+$00 & 7.04e$-$02 & 0.67  & 5.35e$-$04 \\ \cline{2-6}
	&    256 & 1.41e$+$01 & 3.06e$+$00 & 2.76  & 1.64e$-$05 \\ \cline{2-6}
	&   1024 & 5.65e$+$01 & 1.30e$+$01 & 11.4  & 2.24e$-$06 \\ \cline{2-6}
	&   4096 & 2.34e$+$02 & 4.73e$+$01 & 47.0  & 3.10e$-$07 \\ \cline{2-6}
	&  16384 & 1.00e$+$03 & 1.68e$+$02 & 194.0 & 4.36e$-$08 \\
	\hline
\end{tabular}
\caption{\label{tab:DtN_osc} Timing, memory and error results for applying the uniform discetization technique to 
 the locally oscillatory problem with 
 different orders of discretization $n_c$.}
\end{table}
\begin{table}[h!]
\centering
\begin{tabular}{|c|r|c|c|c|c|}
	\hline
	$n_c$ & $N_f$ & $T_{\rm pre}$ & $T_s$ & $R$ & $E_{rel}$ \\
	\hline     
	\multirow{6}{*}{16}
	&      4 & 1.52e$-$01 & 3.35e$-$02 & 0.003 & 1.24e$-$01 \\ \cline{2-6}
	&     16 & 3.73e$-$01 & 1.97e$-$02 & 0.01  & 5.50e$-$03 \\ \cline{2-6}
	&     64 & 9.55e$-$01 & 1.94e$-$02 & 0.06  & 1.08e$-$04 \\ \cline{2-6}
	&    256 & 3.85e$+$00 & 7.24e$-$02 & 0.24  & 1.95e$-$06 \\ \cline{2-6}
	&   1024 & 1.47e$+$01 & 4.40e$-$01 & 1.05  & 1.93e$-$09 \\ \cline{2-6}
	&   4096 & 4.83e$+$01 & 2.21e$+$00 & 4.52  & 1.84e$-$11 \\
    \hline
    \multirow{6}{*}{32}
	&     4 & 8.79e$-$01 & 3.61e$-$02 & 0.04 & 3.20e$-$03 \\ \cline{2-6}                   
	&    16 & 1.12e$+$00 & 2.72e$-$02 & 0.16 & 9.22e$-$05 \\ \cline{2-6}
	&    64 & 3.75e$+$00 & 2.14e$-$01 & 0.67 & 8.69e$-$08 \\ \cline{2-6}
	&   256 & 1.49e$+$01 & 2.46e$+$00 & 2.76 & 2.07e$-$11 \\ \cline{2-6}
	&  1024 & 5.65e$+$01 & 1.33e$+$01 & 11.4 & 1.66e$-$11 \\ \cline{2-6}
	&  4096 & 2.30e$+$02 & 4.78e$+$01 & 47.0 & 6.71e$-$11 \\
	\hline
\end{tabular}
\caption{\label{tab:DtN_wave} Timing, memory and error results for applying the uniform discetization technique to 
 the wave front problem with 
 different orders of discretization $n_c$.}
\end{table}

\subsection{Helmholtz problems}
\label{sec:helm}
This section illustrates the performance of the discretization 
techniques when applied two Helmholtz problems of the form 
\begin{equation*}
    -\Delta u - \omega^2c(\pxx)u = f(x,y)
\end{equation*}
on a square geometry $\Omega $ with an incident wave boundary condition $u(\pxx) = e^{i \omega \vct{d} \cdot \pxx}$ where $\vct{d} = (1,0)$.

Two choices of geometry, coefficient function $c(\pxx) = c(x,y)$ and body load $f(x,y)$ are
considered:\\
\noindent
\textit{Constant coefficient:} For this experiment, $\Omega = (-1,1)^2$ is twenty wavelengths in size ($\omega = 20\pi$),
$c(\pxx) =1$ and $$f(x,y) =  \frac{1}{\sqrt{2\pi\:0.005}}e^{-\frac{{x^2+(y-0.875)^2}}{2 (0.005)^2}}.$$
\noindent
\textit{Variable medium:} For this experiment, $\Omega = (-0.5,0.5)^2$, $\omega = 150$, $f(\pxx) = 0$, and 
$$c(\pxx) = c(x,y) = 4(y-0.2)[1-{\rm erf}(25(|\pxx|-0.3))].$$

\begin{figure}[h!]
 \centering
\setlength{\unitlength}{1mm}
\begin{picture}(95,55)
\put(-10,00){\includegraphics[height = 55mm]{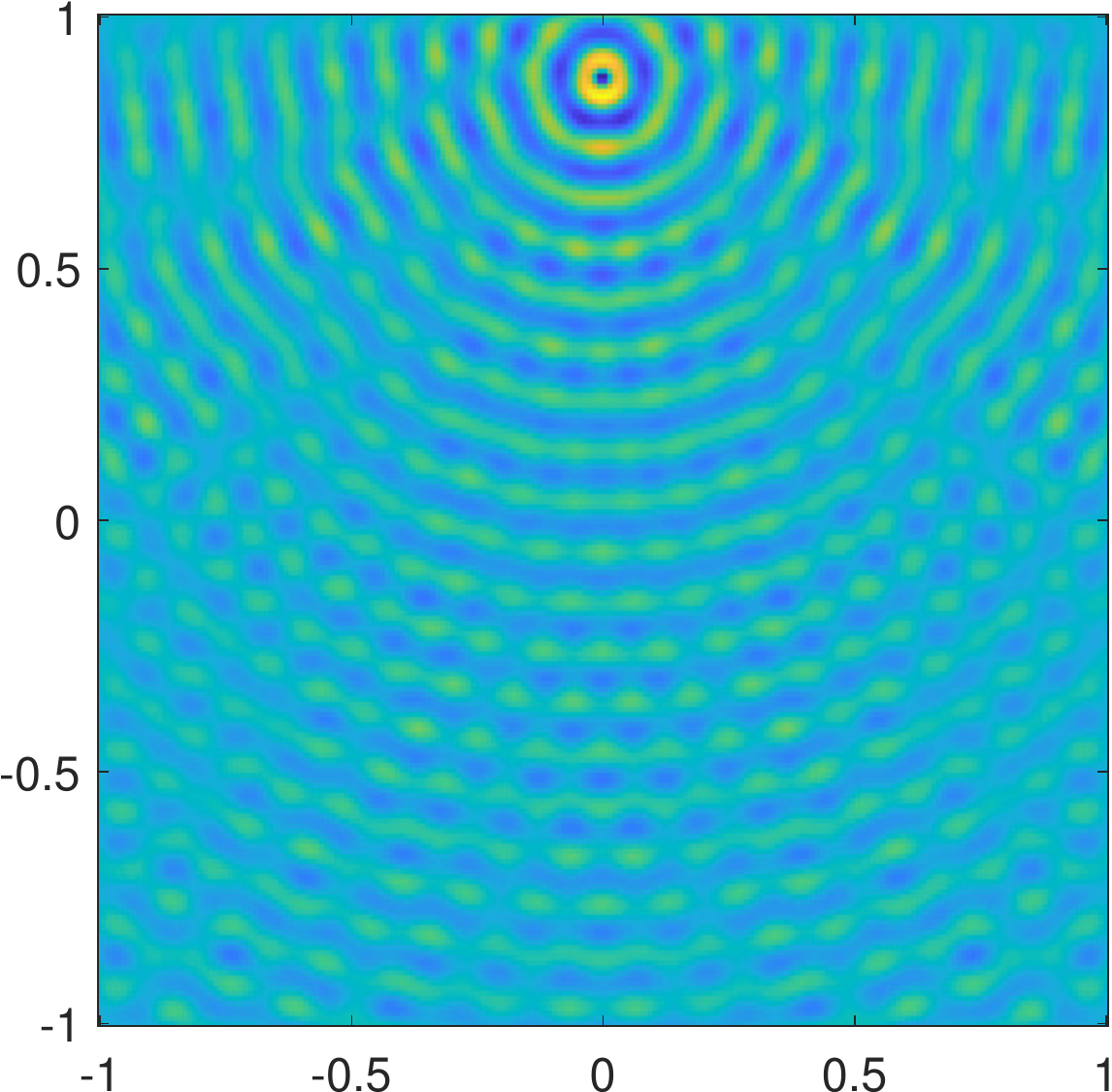}}
\put(21,-3){(a)}
\put(50,00){\includegraphics[height = 55mm]{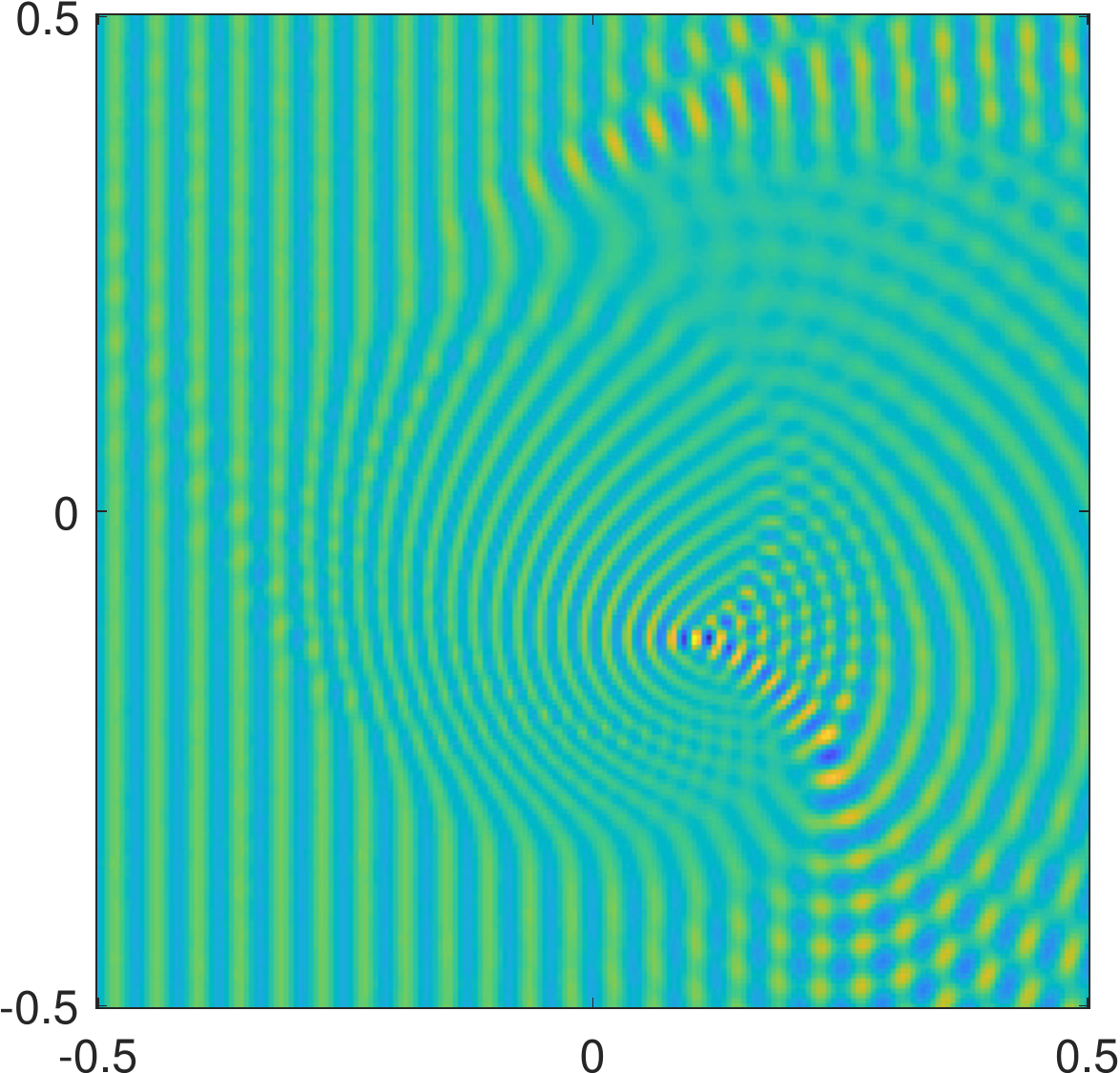}}
\put(81,-3){(b)}
\end{picture}
\caption{\label{fig:ItI_soln} Illustration of the real part of the solutions to the 
problems under consideration in section \ref{sec:helm}: (a) constant coefficient and 
(b) variable coefficient.} 
\end{figure}

Figure \ref{fig:ItI_soln} illustrates the real part of the solutions to these problems.
A reference solution generated by applying the uniform discretization scheme until the 
the relative convergence error was less than the stopping tolerance $\epsilon$ was 
used to generate the reference solution. 

\begin{figure}[h!]
\centering
 \begin{tabular}{ccc}
  \includegraphics[height = 45mm]{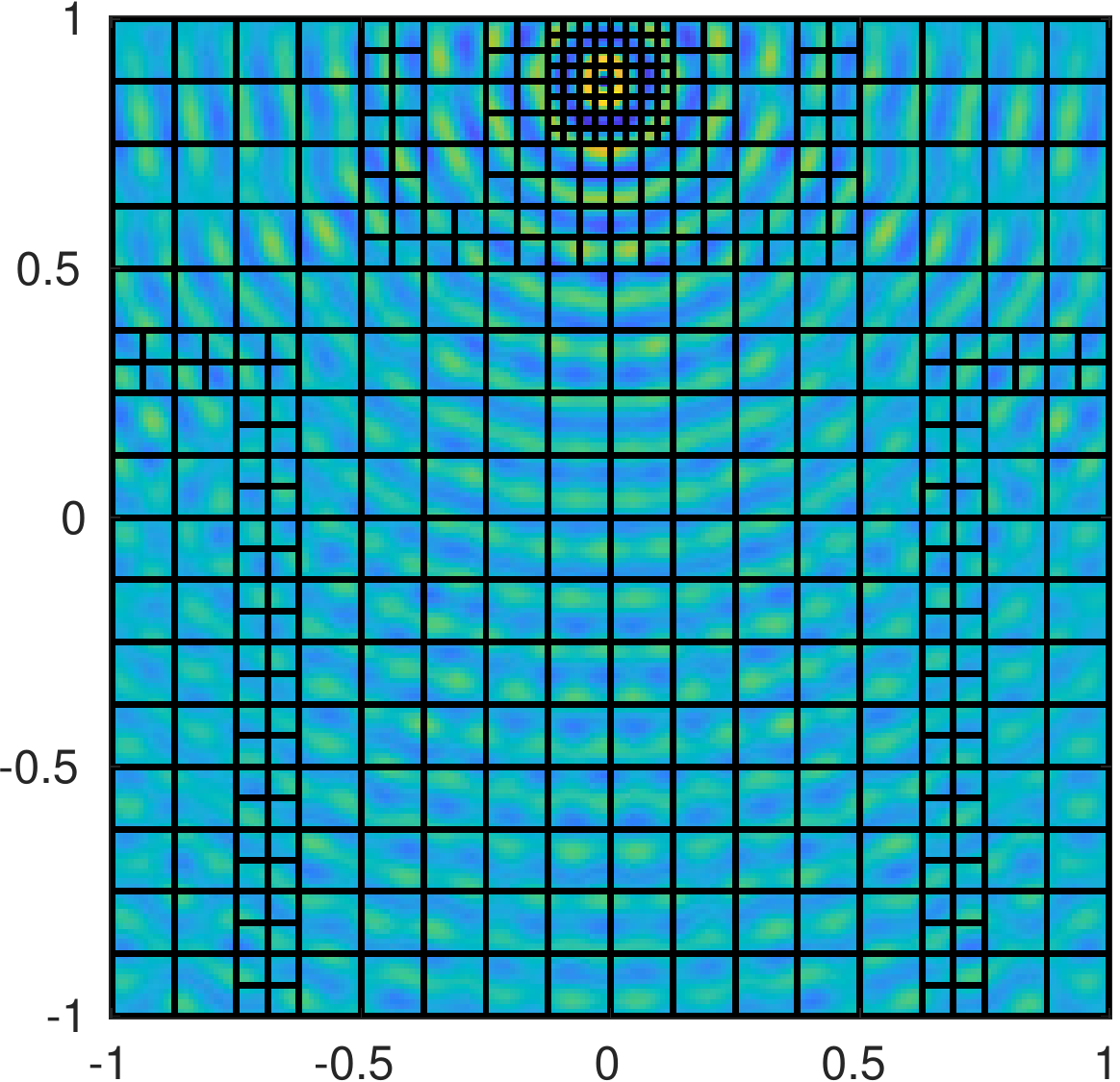} & \mbox{\hspace{20mm}} & 
    \includegraphics[height = 45mm]{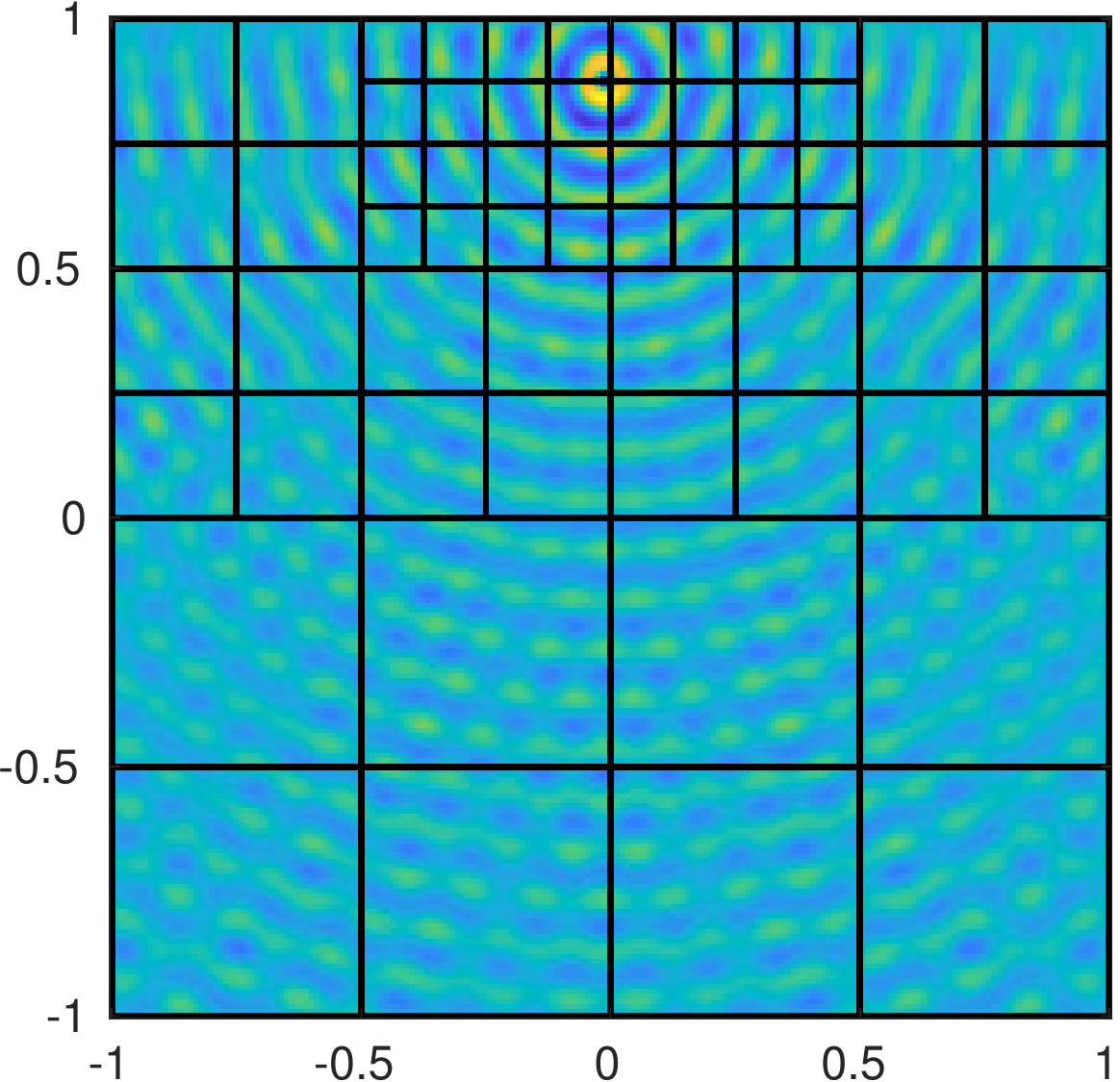}\\
    &&\\
    (a) Constant coefficient with $n_c = 16$ & & (b) Constant coefficient with $n_c = 32$ \\
    \vspace{.2cm} & & \\
  \includegraphics[height = 45mm]{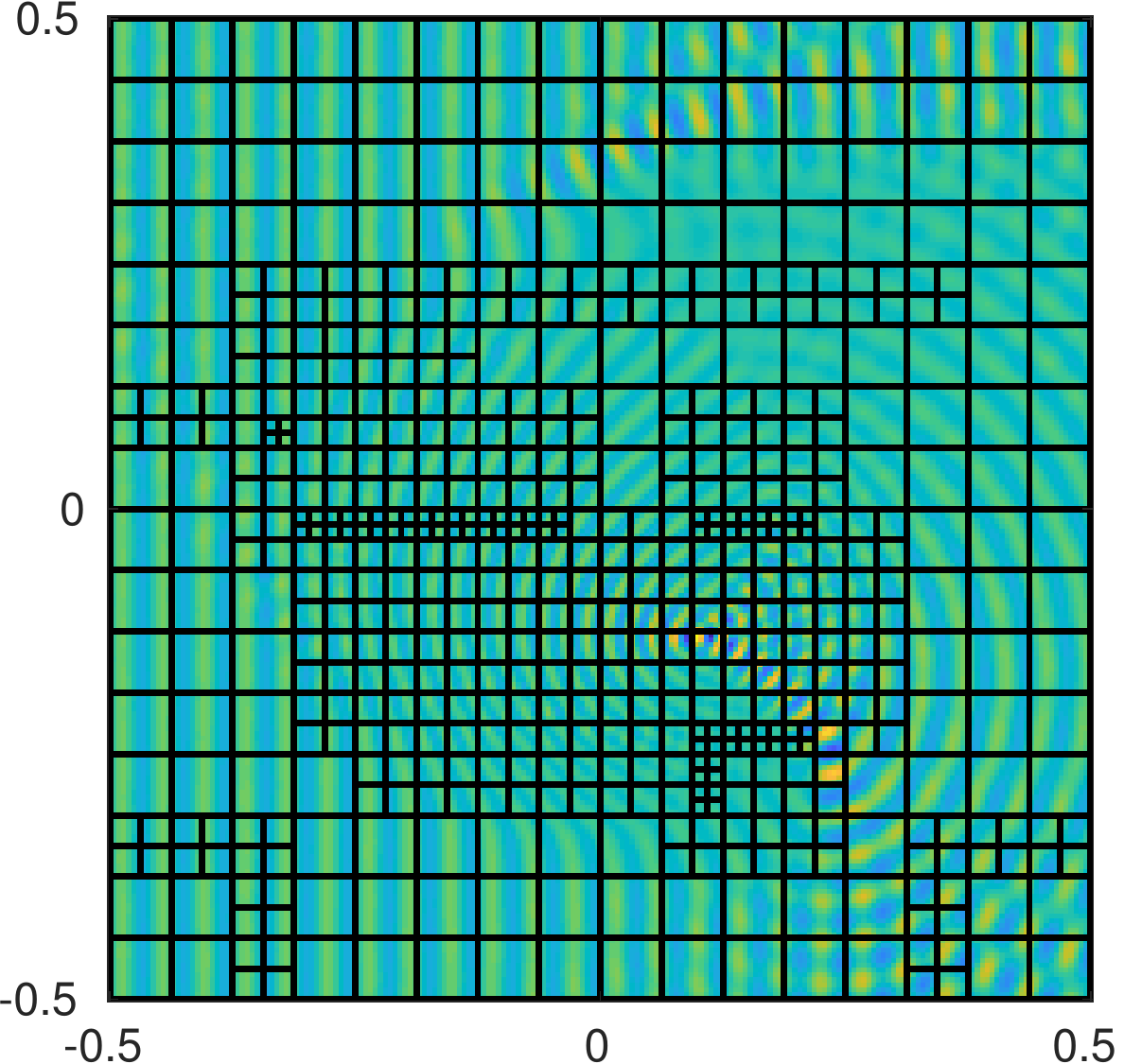} & \mbox{\hspace{20mm}} & 
    \includegraphics[height = 45mm]{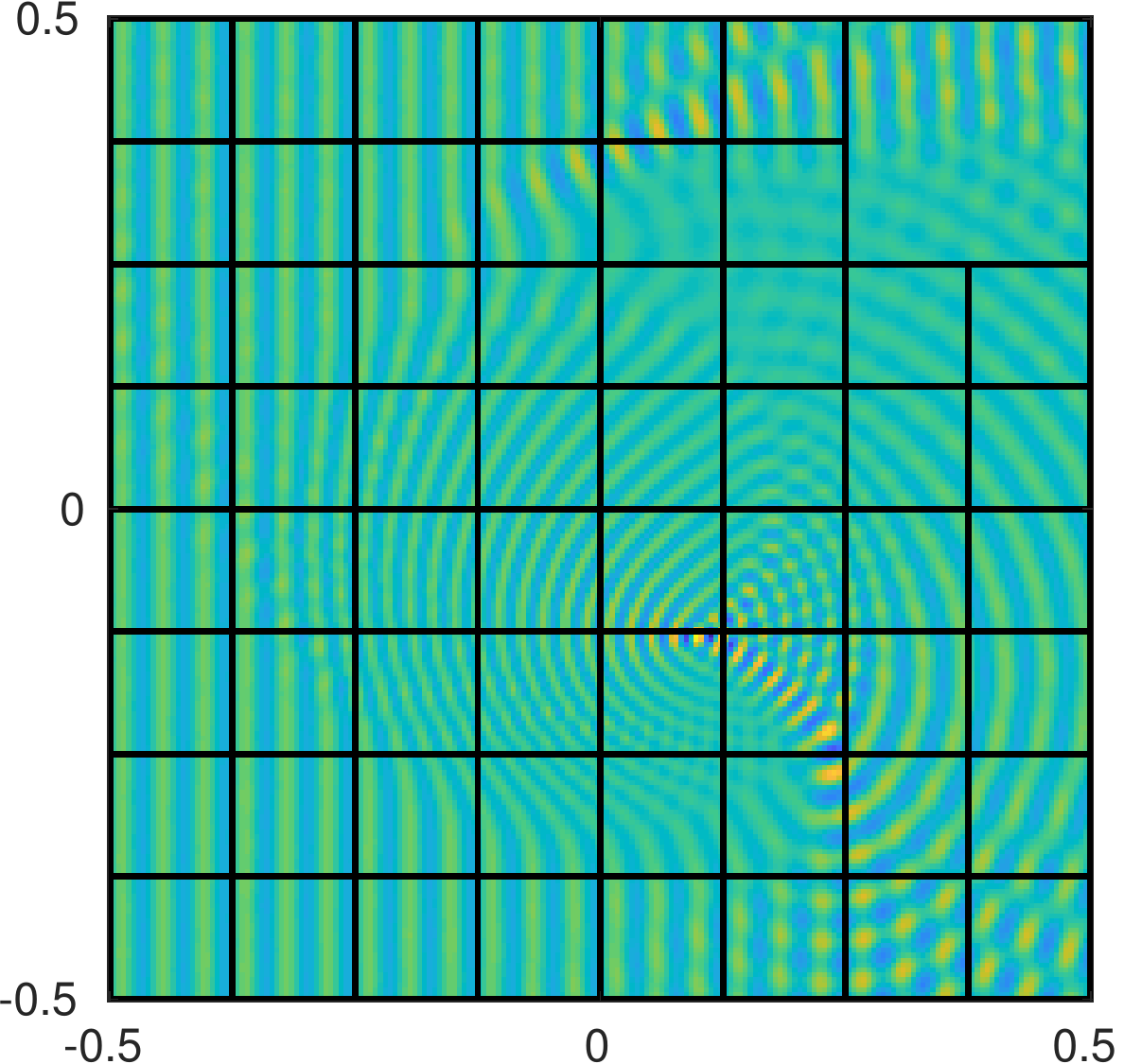}\\
    &&\\
    (c) Variable coefficient with $n_c = 16$ & & (d) Variable coefficient with $n_c = 32$ \\
 \end{tabular}
\caption{\label{fig:helm_grid}  Illustration of the mesh generated by the adaptive discetization
technique (overlayed on the solution) with $n_c =16$ and $n_c = 32$ for the Helmholtz problems
considered in section \ref{sec:helm}.}
\end{figure}

Table \ref{tab:helm_adap} reports on the performance of the adaptive discretization 
technique applied to the Helmholtz problems.  Figure \ref{fig:helm_grid} illustrates
the mesh resulting from the adaptive procedure.  For the constant coefficient case,
the Gaussian body load is not located close to the $n_c=8$ discretization points
on $\Omega^1 = \Omega$, thus the adaptive procedure was not able to capture it.
  To rectify this, we started the adaptive procedure with an 
initialized $16\times 16$ uniform mesh.  
Since high order discretization techniques are better 
suited for high frequency problems, we only consider $n_c = 16$ and $n_c=32$ for the 
variable coefficient problem.  The results for both problems indicate there is 
no benefit in running the adaptive discretization technique with a really high order 
discretization ($n_c = 32$).  The choice of $n_c = 16$ is faster in the precomputation
plus the solve time and memory are comparable.


Tables \ref{tab:const_unif} and \ref{tab:lens_unif} report on the performance of the 
uniform HPS discretization for the constant coefficient and variable coefficient problems,
respectively, with $n_c = 16$ and $32$.    Notice that 
the adaptive method not only requires a smaller number of discretization points, it is
also faster for both the precomputation and apply stages than the uniform method to 
achieve the same accuracy.


\begin{table}[h]
\centering
\begin{tabular}{|c|c|c|c|c|c|c|c|c|}
\hline
Problem& $n_c$ & $N_i$ & $N_f$ & $T_i$ & $T_f$ & $T_s$ & $R$ & $E_{\rm rel}$ \\
 \hline
\multirow{3}{*}{Constant}&$8$  &  1576 & 14470 & 2.54e$+$00 & 3.04e$+$02 & 2.91e$+$00 & 5.242 & 1.02e$-$05 \\ \cline{2-9}
 &$16$ &    1 &   460 & 4.09e$-$02 & 1.71e$+$01 & 2.44e$-$01 & 1.081  & 8.06e$-$06 \\ \cline{2-9}
 &$32$ &    1 &    64 & 1.68e$-$01 & 3.45e$+$01 & 1.31e$-$01 & 1.504  & 9.33e$-$07 \\ \hline	
\multirow{2}{*}{Variable}&$16$ &   64 &    610 & 7.37e$-$01 & 2.06e$+$01 & 3.10e$-$01 & 1.434 & 4.60e$-$06 \\  \cline{2-9}
& $32$ &   16 &     61 & 2.71e$+$00 & 4.31e$+$01 & 1.54e$-$01 & 1.431 & 1.06e$-$05 \\ \hline
 \end{tabular}
\caption{\label{tab:helm_adap} Timing, memory and error results for applying the adaptive technique to the Helmholtz problems in section  \ref{sec:helm}.   }
\end{table}

\begin{table}[h!]
 \centering
\begin{tabular}{|c|r|c|c|c|c|}
	\hline
	$n_c$ & $N_f$ & $T_f$ & $T_s$ & $R$ & $E_{\rm rel}$ \\
	\hline
	\multirow{8}{*}{16}
	&     4 & 1.43e$-$01 & 2.30e$-$02 & 0.007 & 2.50e$-$01 \\ \cline{2-6}
	&    16 & 4.23e$-$01 & 4.25e$-$02 & 0.031 & 7.02e$+$00 \\ \cline{2-6}
	&    64 & 1.31e$+$00 & 3.33e$-$02 & 0.136 & 2.21e$-$01 \\ \cline{2-6}
	&   256 & 5.13e$+$00 & 1.09e$-$01 & 0.591 & 1.40e$-$03 \\ \cline{2-6}
	&  1024 & 2.00e$+$01 & 5.13e$-$01 & 2.554 & 1.20e$-$04 \\ \cline{2-6}
	&  4096 & 7.75e$+$01 & 6.73e$+$00 & 10.98 & 2.02e$-$06 \\ \cline{2-6}
	& 16384 & 3.12e$+$02 & 3.57e$+$01 & 46.98 & 1.91e$-$09 \\ \hline
	\multirow{6}{*}{32}
	&     4 & 1.09e$+$00 & 3.14e$-$02 & 0.087 & 5.35e$-$01 \\ \cline{2-6}
	&    16 & 3.71e$+$00 & 8.50e$-$02 & 0.362 & 1.04e$-$01 \\ \cline{2-6}
	&    64 & 1.40e$+$01 & 3.07e$-$01 & 1.505 & 2.75e$-$02 \\ \cline{2-6}
	&   256 & 5.61e$+$01 & 3.23e$+$00 & 6.238 & 2.20e$-$07 \\ \cline{2-6}
	&  1024 & 2.35e$+$02 & 1.84e$+$01 & 25.83 & 4.05e$-$09 \\ \hline
\end{tabular}
\caption{\label{tab:const_unif} Timing, memory and error results solving the constant 
coefficient Helmholtz problem in section  \ref{sec:helm} with a uniform discretization.   }
\end{table}

\begin{table}[h!]
\centering
\begin{tabular}{|c|r|c|c|c|c|}
	\hline
	$n_c$ & $N_f$ & $T_f$ & $T_s$ & $R$ & $E_{\rm ref}$ \\
	\hline
	\multirow{8}{*}{16}
	&     4 & 1.69e$-$01 & 2.34e$-$02 & 0.007 & 2.35e$-$01 \\ \cline{2-6}
	&    16 & 5.00e$-$01 & 1.58e$-$02 & 0.031 & 2.89e$-$01 \\ \cline{2-6}
	&    64 & 1.49e$+$00 & 2.63e$-$02 & 0.136 & 1.89e$-$01 \\ \cline{2-6}
	&   256 & 5.91e$+$00 & 1.01e$-$01 & 0.591 & 7.76e$-$04 \\ \cline{2-6}
	&  1024 & 2.43e$+$01 & 8.95e$-$01 & 2.554 & 7.47e$-$08 \\ \cline{2-6}
	&  4096 & 9.04e$+$01 & 8.98e$+$00 & 10.98 & 1.24e$-$10 \\ \cline{2-6}
	& 16384 & 3.21e$+$02 & 3.77e$+$01 & 46.98 & 7.78e$-$10 \\ \hline
	\multirow{6}{*}{32}
	&     4 & 1.21e$+$00 & 3.20e$-$02 & 0.087 & 2.64e$-$01 \\ \cline{2-6}
	&    16 & 4.45e$+$00 & 3.67e$-$02 & 0.362 & 1.54e$-$01 \\ \cline{2-6}
	&    64 & 1.63e$+$01 & 2.85e$-$01 & 1.505 & 2.29e$-$06 \\ \cline{2-6}
	&   256 & 6.19e$+$01 & 3.72e$+$00 & 6.238 & 9.24e$-$11 \\ \cline{2-6}
	&  1024 & 2.76e$+$02 & 1.71e$+$01 & 25.83 & 6.23e$-$10 \\ \hline
\end{tabular}
\caption{\label{tab:lens_unif} Timing, memory and error results solving the variable 
coefficient Helmholtz problem in section  \ref{sec:helm} with a uniform discretization.   }
\end{table}

\section{Concluding remarks}
\label{sec:summary}
This manuscript presents an high order adaptive discretization technique that 
comes with an efficient direct solver.  The HPS method presented
here uses a \textit{new} local pseudospectral discretization that does not involve 
corner points.  By removing the corner points, the leaf computations are less 
expensive and more stable than the previous version of the method.  

The adaptive discretization technique utilizes the modified local tensor product 
basis to look at convergence of the directional Chebychev expansions to 
determine which regions of the domain $\Omega$ need refinement.  Since a discretization
is based on decomposing the domain, updating the accompanying direct solver 
after refinement is inexpensive.  The numerical results show that method is 
able to achieve the user prescribed accuracy and refines only of regions where 
it is necessary. For all problems the cost of applying and storing the direct solver
resulting from the adaptive discretization technique is less than using a uniform 
discretization.  For problems where the solution is globally oscillatory the 
cost of adaptive discretization technique is less than a uniform discretization.

\section{Acknowledgements} The work by A. Gillman is supported by the Alfred P. Sloan foundation and the
National Science Foundation (DMS-1522631). 
The work by P. Geldermans is supported by the National Science Foundation Graduate 
Research Fellowship under Grant No. 1450681.


\section{Appendix}
\label{sec:app}
Consider the variable coefficient Helmholtz impedance boundary value problem 
\begin{equation}\begin{aligned}
    -\Delta u - \omega^2c(\pxx)u &= s(x,y) \qquad & \pxx\in\Omega\\
 \frac{\partial u}{\partial \nu} +i\eta u &= t(x,y) \qquad & \pxx\in \partial\Omega = \Gamma.
 \end{aligned}
 \label{eq:varyHelm}
 \end{equation}
where $\nu$ denotes the outward facing normal vector, $c(\pxx)$ is a smooth 
function,  $\omega\in\mathbb{R}$, and $\eta\in\mathbb{C}$. 

This section presents the technique for solving variable coefficient Helmholtz problems 
such as (\ref{eq:varyHelm}) using the HPS method. This technique uses impedance-to-impedance
(ItI) operators instead of the DtN operators used in the body of the paper.  

\begin{definition}[impedance-to-impedance map]
Fix $\eta\in\mathbb{C}$, and $\mathcal{Re} \eta\neq 0$. Let
\begin{eqnarray}
f&:=& u_n+i\eta u|_\Gamma
\label{f}
\\
g&:=&u_n-i\eta u|_\Gamma
\label{g}
\end{eqnarray}
be Robin traces of $u$. We refer to $f$ and $g$ as the ``incoming'' and ``outgoing'' (respectively)
impedance data.
For any $\omega>0$, the \textit{ItI operator} $R:L^2(\Gamma)\to L^2(\Gamma)$ is defined by
\begin{equation}
R f = g
\label{R}
\end{equation}
for $f$ and $g$ the Robin traces of $u$ the solution of (\ref{eq:varyHelm}),
for all $f \in L^2(\Gamma)$.
 
\end{definition}

To make the solution technique useful for different choices of $s(x,y)$, we choose
to represent the solution $u$ as a superposition of the homogeneous solution $w$
and the particular solution $z$; i.e. $u = w+z$ where $z$ is the solution 
of the following boundary value problem
\begin{align*}
    -\Delta z - \omega^2c(\pxx)z &= s(x,y) \qquad & \pxx\in\Omega\\
 \frac{\partial z}{\partial \nu} +i\eta z &= 0 \qquad & \pxx\in \partial\Omega = \Gamma
\end{align*}
and $w$ is the solution of 
\begin{align*}
    -\Delta w - \omega^2c(\pxx)w &= 0 \qquad & \pxx\in\Omega\\
 \frac{\partial w}{\partial \nu} +i\eta w &= t(x,y) \qquad & \pxx\in \partial\Omega = \Gamma.
\end{align*}

Section \ref{sec:ItIleaf} presents the leaf computation
and section \ref{sec:ItImerge} presents the technique for merging two boxes.  Throughout 
the notation is kept consistent with that of section \ref{sec:HPS}. 
When there is no body load (i.e. $s(x,y)= 0$), the method from \cite{2013_martinsson_ItI} is recovered.

\subsection{Leaf computation}
\label{sec:ItIleaf}
This section presents the construction of the homogeneous and particular solutions to (\ref{eq:varyHelm})
using the modified spectral collocation method from section \ref{sec:leaf}.  Additionally, 
a matrix $\mtx{R}$ approximating the ItI operator for the homogeneous boundary value problem and
the impedance boundary data generated by the particular solution are constructed.

Let $\mtx{N}$ denote the matrix that takes normal derivatives of the basis functions.  Then 
$\mtx{N}$ is given by 

$$\mtx{N} = \left[\begin{array}{c} -\mtx{D}_x(I_s,I^\tau) \\
                           \mtx{D}_y(I_e,I^\tau)\\
                           \mtx{D}_x(I_n,I^\tau)\\
                           -\mtx{D}_y(I_w,I^\tau)\end{array}\right].$$
                           
Then the matrix for creating the incoming impedance data is 
$$\mtx{F} = \mtx{N}+i\eta\mtx{I}_{n_c^2}(I_b,I^\tau)$$
and the matrix for creating the outgoing impedance data is 
$$\mtx{G} = \mtx{N}-i\eta\mtx{I}_{n_c^2}(I_b,I^\tau)$$
where $\mtx{I}_{n_c^2}$ is the identity matrix of size 
$n_c^2$.

Then the discretized body load problem to find the approximation to $z$ at the 
collocation points takes the form

\begin{equation}
 \mtx{B}\vtwo{\vct{z}_{b}}{\vct{z}_i} = \vtwo{\mtx{F}}{\mtx{A}(i,b) \ \mtx{A}(i,i)} \vct{z} = \vtwo{\mtx{0}}{\vct{s}}
 \label{eq:particular}
 \end{equation}
where $\vct{z}$ is the vector with the approximate values of $z$ at the collocation points,
and $\vct{s}$ is $s(x,y)$ evaluated at the interior points. 

So the solution operator $\mtx{Y}$ which gives the approximate particular solution is the 
solution to 
$$\mtx{B}\mtx{Y} = \vtwo{\vct{0}_{4n_c-4\times (n_c-2)^2}}{\vct{I}_{(n_c-2)^2}}.$$

Likewise the solution operator $\mtx{\Psi}$ which give the approximate solution to the 
homogeneous problem is the solution to 
$$\mtx{B}\mtx{\Psi} = \vtwo{\vct{I}_{4n_c-4}}{\vct{0}_{ (n_c-2)^2\times 4n_c-4}}.$$

To construct the outgoing impedance data from the particular solution $\mtx{h}$, 
the matrix $\mtx{G}$ needs to be applied to the solution of (\ref{eq:particular});
i.e.
$$\vct{h} = \mtx{G}\mtx{Y}\vtwo{\mtx{0}}{\vct{s}} = \mtx{W}\vtwo{\mtx{0}}{\vct{s}}.$$

The approximate ItI operator is constructed in the same manner as in \cite{2013_martinsson_ItI}.  That 
is 

$$\mtx{R} = \mtx{G}\mtx{\Psi}.$$

Putting these together, we find that the outgoing impedance data from the box is 
given by 
$$\vct{g} = \mtx{R}\vct{t} + \vct{h}$$
where $\vct{t}$ is the evaluation of the incoming boundary data function $t(x,y)$ at 
the points on the boundary.

\subsection{Merge two boxes}
\label{sec:ItImerge}
This section presents the technique for merging two boxes $\Omega^\tau = \Omega^\alpha\cup\Omega^\beta$ 
for which the ItI matrices and outgoing impedance data from the particular 
solution has already been computed.  In other words, the matrices $\mtx{R}^\alpha$ and $\mtx{R}^\beta$ 
along with the vectors $\vct{h}^\alpha$ and $\vct{h}^\beta$ are available.  For consistency, we used
the same notation as in \cite{2013_martinsson_ItI}.  In this section, it is important to 
remember that the unlike the DtN version of the algorithm, the normal derivatives are 
always pointing exterior to the region they are defined on.

Using the same ordering as in section \ref{sec:merge}, the outgoing impedance data 
for boxes $\alpha$ and $\beta$ are given by

$$  \vtwo{\vct{g}^\alpha_1}{\vct{g}^\alpha_3}=\mtwo{\mtx{R}^\alpha_{11}}{\mtx{R}^\alpha_{13}}{\mtx{R}^\alpha_{31}}{\mtx{R}^\alpha_{33}} 
\vtwo{\vct{t}^\alpha_1}{\vct{t}^\alpha_3} 
+ \vtwo{\vct{h}^\alpha_1}{\vct{h}^\alpha_3} ; \qquad
\vtwo{\vct{g}^\beta_2}{\vct{g}^\beta_3} = 
\mtwo{\mtx{R}^\beta_{22}}{\mtx{R}^\beta_{23}}{\mtx{R}^\beta_{32}}{\mtx{R}^\beta_{33}} 
\vtwo{\vct{t}^\beta_2}{\vct{t}^\beta_3}
 + \vtwo{\vct{h}^\beta_2}{\vct{h}^\beta_3}$$

where $\vtwo{\vct{h}^\alpha_1}{\vct{h}^\alpha_3}$ and $\vtwo{\vct{h}^\beta_2}{\vct{h}^\beta_3}$
are the outgoing impedance data due to the particular solutions on each box.

Since the normal vectors are opposite in each box, we know $\vct{t}^\alpha_3 = -\vct{g}^\beta_3$ and $\vct{g}^\alpha_3 = -\vct{t}^\beta_3$.
Using this information in the bottom row equations, $\vct{t}^\alpha_3$ and $\vct{t}^\beta_3$ can 
found in terms of $\vct{t}_1^\alpha$, $\vct{t}_2^\beta$, $\vct{h}^3_\alpha$, and $\vct{h}^3_\beta$.
They are given by
\begin{equation}
\vct{t}_3^\alpha = \mtx{W}^{-1}\left[\mtx{R}_{33}^\beta \mtx{R}_{31}^\alpha | -\mtx{R}_{32}^\beta\right] 
\vtwo{\vct{t}_1^\alpha}{\vct{t}_{2}^\beta} +\mtx{W}^{-1}\left(\mtx{R}_{33}^\beta\vct{h}^\alpha_3-\vct{h}_3^\beta\right)
 \label{eq:talpha}
\end{equation}
and
\begin{equation}
\vct{t}^\beta_3 = \left[-\mtx{R}_{31}^\alpha -\mtx{R}^\alpha_{33}\mtx{W}^{-1}\mtx{R}_{33}^\beta\mtx{R}^\alpha_{31} | \mtx{R}_{33}^\alpha \mtx{W}^{-1}\mtx{R}_{32}^\beta\right]
\vtwo{\vct{t}_{1}^\alpha}{\vct{t}_{2}^\beta} -\left(\mtx{I}+\mtx{R}^\alpha_{33}\mtx{W}^{-1}\mtx{R}_{33}^\beta\right)\vct{h}^\alpha_3 +\mtx{R}_{33}^\alpha \mtx{W}^{-1}\vct{h}^\beta_3
 \label{eq:tbeta}
\end{equation}
where $\mtx{W} = \mtx{I}-\mtx{R}_{33}^\beta\mtx{R}^\alpha_{33}$.  

Plugging (\ref{eq:talpha}) and (\ref{eq:tbeta}) into the top row equations results in the following expression 
for the outgoing impedance data for the box $\Omega^\alpha\cup \Omega^\beta$

\begin{equation}\begin{split}
 \vtwo{\vct{g}^\alpha_1}{\vct{g}^\beta_2}& = 
 \mtwo{\mtx{R}_{11}^\alpha +\mtx{R}^\alpha_{13}\mtx{W}^{-1}\mtx{R}_{33}^\beta\mtx{R}_{31}^\alpha}{-\mtx{R}_{13}^\alpha \mtx{W}^{-1}\mtx{R}_{32}^\beta}
 {-\mtx{R}_{23}^\beta\left(\mtx{R}^\alpha_{31}+\mtx{R}_{33}^\alpha\mtx{W}^{-1}\mtx{R}_{33}^\beta\mtx{R}_{31}^\alpha\right)}
 {\mtx{R}_{22}^\beta+\mtx{R}_{23}^\beta\mtx{R}_{33}^\alpha\mtx{W}^{-1}\mtx{R}_{32}^\beta}
 \vtwo{\vct{t}_1^\alpha}{\vct{t}_2^\beta} \\
 & \qquad +\vtwo{\vct{h}^\alpha_1}{\vct{h}_2^\beta} +
 \vtwo{\mtx{R}_{13}^\alpha \mtx{W}^{-1}\left(\mtx{R}_{33}^\alpha\vct{h}^\alpha_3-\vct{h}^\beta_3\right)}
 {-\mtx{R}^\beta_{23}\left(\mtx{I}+\mtx{R}_{33}^\alpha\mtx{W}^{-1}\mtx{R}_{33}^\beta\right)\vct{h}^\alpha_3+\mtx{R}_{23}^\beta\mtx{R}_{33}^\alpha\mtx{W}^{-1}\vct{h}^\beta_3}\\
 &= \mtx{R}^\tau  \vtwo{\vct{t}_1^\alpha}{\vct{t}_2^\beta} 
 +\vtwo{\hat{\vct{h}}^\alpha_1}{\hat{\vct{h}}_2^\beta} 
 \end{split}
 \label{eq:out}
\end{equation}

\subsection{The full algorithm}

As with the homogeneous DtN solution technique, the solver can be broken into the 
precomputation and the solve phase.  The precomputation for a leaf box $\tau$ is 
similar to before except now a solution operator $\mtx{Y}^\tau$ yielding the particular solution
on $\tau$ and a matrix $\mtx{W}$ giving the outgoing particular impedance data are constructed.
Also, instead of a DtN matrix, an ItI matrix is constructed.  The precomputation for
a box $\tau$ with children $\alpha$ and $\beta$ is more intense.  A collection of operators
giving the incoming impedance data on the shared edge are constructed from the 
incoming impedance data from $\tau$ and the outgoing impedance particular solution data (which is 
not yet computed) on that edge from both $\alpha$ and $\beta$.  Thus a collection of 
operators for constructing the outgoing particular solution on the shared edge are 
constructed as well as the operators needed to construct the outgoing impedance particular solution 
data on the boundary of $\tau$.  Notice looking the formulas (\ref{eq:talpha}), (\ref{eq:tbeta}) and 
(\ref{eq:out}) there is significant overlap in computation thus keeping the cost and 
memory of the precomputation in check.

The solve step sweeps the tree twice (instead of once as in the homogeneous solver).  
First, starting from the leaf boxes moving up the tree to $\Omega^1$, the outgoing 
impedance particular solution data are constructed.  Then using this information 
along with the boundary condition on $\Omega$, the incoming impedance boundary 
data is propagated from the top of the tree down to the leaf boxes.

\end{document}